%% file: thesis.tex
\documentclass[titlepage,oneside,11pt]{report}

\usepackage{amsmath}
\usepackage{amsfonts}
\usepackage{vmargin}
\usepackage{amsthm}
\usepackage{amscd}
\usepackage{epsfig}

\include{defs}

\setpapersize{USletter}
\normalsize 
\setmarginsrb{1.75in}{1in}{1.25in}{1in}
	{\baselineskip}{18pt}{\baselineskip}{18pt}

\title{Integral Expressions for the Vassiliev Knot Invariants}
\author{Dylan Thurston}
\date{April 11, 1995\\(Figures redrawn January 28, 1999)}

\begin{document}
\bibliographystyle{plain}
\maketitle


\input{abstract}

\tableofcontents
\input{intro}

\input{vassiliev}
\input{configspace}

\input{tying}
\input{nastyface}
\input{anomaly}
\appendix
\input{linalg}
\input{graph}
\nocite{*}
\bibliography{knots}
\end{document}

%% file: defs.tex
\newcommand\RR{\mathbb{R}}
\newcommand\Vect{\text{Vect}}
\newcommand\Bl{\mathrm{Bl}}
\newcommand\oC{C^0}
\newcommand\calT{\mathcal{T}}

\newcommand\fS{\mathfrak{S}}

\newcommand\orient{\operatorname{or}}
\newcommand\supp{\operatorname{supp}}
\newcommand\sgn{\operatorname{sgn}}

\DeclareMathAlphabet\EuScript{U}{eus}{m}{n}
\SetMathAlphabet\EuScript{bold}{U}{eus}{b}{n}

\newcommand\knots{\EuScript{K}}
\newcommand\fchain{\EuScript{F}}

\usepackage{calc}
\def\mathcenter#1{
	\raisebox{.5ex+\depth-.5\totalheight}{\hbox{#1}}
}
\newcommand\intersect{\mathcenter{\input{draws/intersect.tex}}}
\newcommand\overcross{\mathcenter{\input{draws/overcross.tex}}}
\newcommand\undercross{\mathcenter{\input{draws/undercross.tex}}}

\newcommand\A{{K_1}}
\newcommand\B{{K_2}}

\newtheorem*{theorem}{Theorem}
\newtheorem{proposition}{Proposition}[chapter]
\newtheorem{lemma}[proposition]{Lemma}

\theoremstyle{definition}
\newtheorem{definition}[proposition]{Definition}

\newtheorem{case}{Case}

%% file: draws/intersect.tex
\begin{picture}(0,0)%
\epsfig{file=draws/intersect.pstex}%
\end{picture}%
\setlength{\unitlength}{3947sp}%
\begingroup\makeatletter\ifx\SetFigFont\undefined
\def\x#1#2#3#4#5#6#7\relax{\def\x{#1#2#3#4#5#6}}%
\expandafter\x\fmtname xxxxxx\relax \def\y{splain}%
\ifx\x\y   
\gdef\SetFigFont#1#2#3{%
  \ifnum #1<17\tiny\else \ifnum #1<20\small\else
  \ifnum #1<24\normalsize\else \ifnum #1<29\large\else
  \ifnum #1<34\Large\else \ifnum #1<41\LARGE\else
     \huge\fi\fi\fi\fi\fi\fi
  \csname #3\endcsname}%
\else
\gdef\SetFigFont#1#2#3{\begingroup
  \count@#1\relax \ifnum 25<\count@\count@25\fi
  \def\x{\endgroup\@setsize\SetFigFont{#2pt}}%
  \expandafter\x
    \csname \romannumeral\the\count@ pt\expandafter\endcsname
    \csname @\romannumeral\the\count@ pt\endcsname
  \csname #3\endcsname}%
\fi
\fi\endgroup
\begin{picture}(324,324)(1789,-1273)
\end{picture}

%% file: draws/overcross.tex
\begin{picture}(0,0)%
\epsfig{file=draws/overcross.pstex}%
\end{picture}%
\setlength{\unitlength}{3947sp}%
\begingroup\makeatletter\ifx\SetFigFont\undefined
\def\x#1#2#3#4#5#6#7\relax{\def\x{#1#2#3#4#5#6}}%
\expandafter\x\fmtname xxxxxx\relax \def\y{splain}%
\ifx\x\y   
\gdef\SetFigFont#1#2#3{%
  \ifnum #1<17\tiny\else \ifnum #1<20\small\else
  \ifnum #1<24\normalsize\else \ifnum #1<29\large\else
  \ifnum #1<34\Large\else \ifnum #1<41\LARGE\else
     \huge\fi\fi\fi\fi\fi\fi
  \csname #3\endcsname}%
\else
\gdef\SetFigFont#1#2#3{\begingroup
  \count@#1\relax \ifnum 25<\count@\count@25\fi
  \def\x{\endgroup\@setsize\SetFigFont{#2pt}}%
  \expandafter\x
    \csname \romannumeral\the\count@ pt\expandafter\endcsname
    \csname @\romannumeral\the\count@ pt\endcsname
  \csname #3\endcsname}%
\fi
\fi\endgroup
\begin{picture}(324,331)(2089,-1280)
\end{picture}

%% file: draws/undercross.tex
\begin{picture}(0,0)%
\epsfig{file=draws/undercross.pstex}%
\end{picture}%
\setlength{\unitlength}{3947sp}%
\begingroup\makeatletter\ifx\SetFigFont\undefined
\def\x#1#2#3#4#5#6#7\relax{\def\x{#1#2#3#4#5#6}}%
\expandafter\x\fmtname xxxxxx\relax \def\y{splain}%
\ifx\x\y   
\gdef\SetFigFont#1#2#3{%
  \ifnum #1<17\tiny\else \ifnum #1<20\small\else
  \ifnum #1<24\normalsize\else \ifnum #1<29\large\else
  \ifnum #1<34\Large\else \ifnum #1<41\LARGE\else
     \huge\fi\fi\fi\fi\fi\fi
  \csname #3\endcsname}%
\else
\gdef\SetFigFont#1#2#3{\begingroup
  \count@#1\relax \ifnum 25<\count@\count@25\fi
  \def\x{\endgroup\@setsize\SetFigFont{#2pt}}%
  \expandafter\x
    \csname \romannumeral\the\count@ pt\expandafter\endcsname
    \csname @\romannumeral\the\count@ pt\endcsname
  \csname #3\endcsname}%
\fi
\fi\endgroup
\begin{picture}(324,331)(1789,-1280)
\end{picture}

%% file: abstract.tex
\begin{abstract}
It has been folklore for several years in the knot theory community
that certain integrals on configuration space, originally motivated by
perturbation theory
for the Chern-Simons field theory, converge and yield knot invariants.  This
was proposed independently by Gaudagnini, Martellini, and
Mintchev~\cite{gmm:PertChern} and Bar-Natan~\cite{bn:phd}.  The
analytic difficulties involved in proving convergence and invariance
were reportedly worked out by Bar-Natan~\cite{bn:phd,bn:pcs},
Kontsevich~\cite{kontsevich:vassiliev,kontsevich:graphs},
and Axelrod and
Singer~\cite{Axelrod-Singer:ChernPert,Axelrod-Singer:ChernPertII}.

But I know of no elementary exposition of this fact.
Bar-Natan~\cite{bn:pcs} only proves invariance for the degree 2 invariant.
Kontsevich's
exposition is decidedly non-elementary and leaves many details
implicit.  Axelrod and Singer's papers take a physics point of
view (and are thus difficult for mathematicians to read)
and only discuss the related invariants for 3-manifolds.
More recently, Bott and Taubes~\cite{bt:knot} have explained the degree 2
invariant from a purely topological point of view, but again omit the
higher degree cases.

This thesis is an attempt to remedy this
lack.
I adopt an almost exclusively topological point of view, rarely
mentioning Chern-Simons theory.  For an explanation of the physics
involved, see, for instance, Bar-Natan~\cite{bn:pcs}.

There are also a few new results in this thesis.  These include a new
construction of the functorial compactification of configuration space
(Section~\ref{sec:configdef}) as well as
some variations on the integrals.  For a
suitable choice of this variation, the integral reduces to counting
\emph{tinkertoy diagrams} (Section~\ref{sec:tinkertoy}).  In
particular, the invariants constructed
take values in $\mathbb{Q}$.
\end{abstract}

%% file: intro.tex
\chapter{Introduction} \label{sec:introduction}
\section{Motivation}
\label{sec:motivation}
First, I'll present some motivation for introducing these invariants,
primarily intended for topologists.

A (mathematical) {\em knot}%
\footnote{As opposed to a physical knot, which has a finite thickness
and has ends.}
is a loop in $\RR^3$:
an embedding $K: S^1 \rightarrow \RR^3$.  An embedding $K$ is a smooth
map with two conditions:
\begin{enumerate}
\item Immersion condition:
$\frac{dK}{d\theta}(\theta) \ne 0$, $\forall \theta \in S^1$
\item Non self-intersecting:
$K(\theta_0) \ne K(\theta_1)$, $\forall \theta_0,\theta_1 \in S^1$
\end{enumerate}

The first condition is purely local and not terribly interesting.  (In
fact, removing that condition yields an essentially
equivalent space.)  The
second condition is much more interesting.  It's a global condition and is
what makes the study of knots much harder than the study of loops.
How can we reduce it to a finite condition?  One answer comes from the
tradition of raising maps to the product of the spaces involved.
Observe that from the map $K$, we can construct
\begin{equation}
K^{\times n}: (S^1)^n \rightarrow (\RR^3)^n \cr
(\theta_1, \ldots, \theta_n) \mapsto
	(K(\theta_1), \ldots, K(\theta_n))
\end{equation}
This may seem like a trivial remark, but it has far-reaching consequences.  In
particular, because of condition 2, $K^{\times n}$ restricts to a map between
the configuration spaces.  The {\em configuration space} $\oC_n(M)$ of
$n$ points 
in a manifold $M$ is the space of $n$ {\em distinct} points
in $M$.  Formally,
\begin{equation}
\oC_n(M) = \{ (x_1, x_2, \ldots, x_n) \in M^n: x_i \ne x_j
	\mathrm{\ for\ } i\ne j\}
\end{equation}
and, if $f: M \rightarrow N$ is an embedding, we can define
\begin{equation}
\oC_n(f): \oC_n(M) \rightarrow \oC_n(N)
\end{equation}
as the restriction of $f^{\times n}$ to $\oC_n(M)$.
You can therefore try to gain information about the map $f$ by looking
at the map on the configuration spaces $f:C^n(M) \rightarrow C^n(N)$.

In this thesis, I begin to carry out this program in the case of knots.
$\oC_n(\RR^3)$ has a useful set of homology generators given by the
pullbacks of a 2-form on $S^2$ via the direction
map for each pair of points in $\oC_n(\RR^3)$.  By pulling back these
forms to $\oC_n(S^1)$ and taking sums so that the boundary terms
cancel, I will construct explicit formulas for the Vassiliev
knot invariants.

\section{Linking Number}

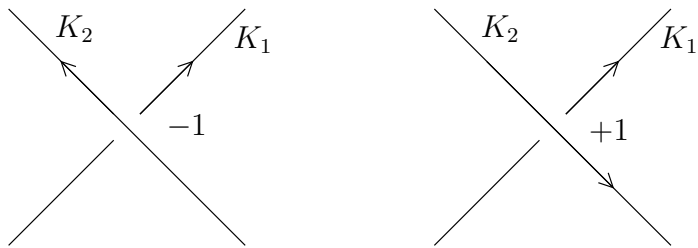
\begin{figure}
\centerline{\input{draws/link-sign.tex}}
\caption{Sign convention for the linking number.  The link is viewed
from directly above.}\label{fig:crossing}
\end{figure}

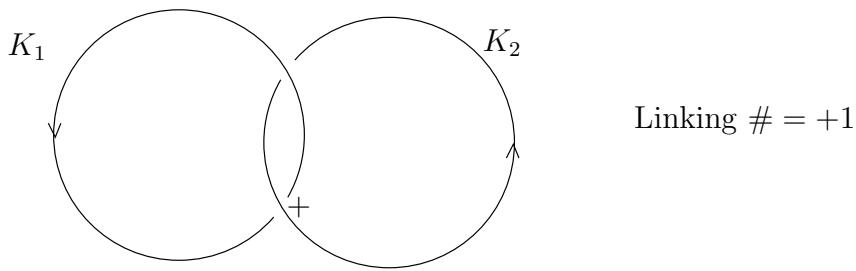
\begin{figure}
\centerline{\input{draws/link-Hopf.tex}}
\caption{Computation of the linking number for
the Hopf link}\label{fig:linksample}
\end{figure}

To help explain these ideas, I'll give a simple example.
The \emph{linking number} is a well-known invariant of two-component
links.  (An $n$-component \emph{link} is an embedding of a disjoint
union of $n$ circles in $\RR^3$.)  Label the two circles $\A$
and $\B$ and consider the set of \emph{$\A\B$-crossings}---points of $\B$
that lie directly above
a point of $\A$.
  If the link
is in general position, there will only be a finite number of
crossings, and they will project to distinct points in the $xy$ plane.
Count up the crossings with signs, as in Figure~\ref{fig:crossing}.
A crossing
is counted with a positive sign if $\B$ runs from left to right as you
look along the direction of $\A$, negative otherwise.  See
Figure~\ref{fig:linksample} for an example of this computation of the
linking number.

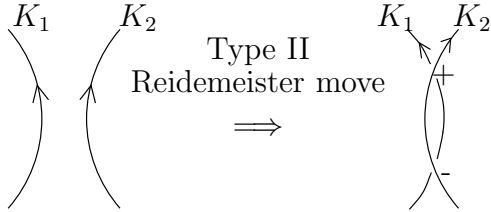
\begin{figure}
\centerline{\input{draws/link-inv.tex}}
\caption{The linking number is an invariant.}\label{fig:linkinvar}
\end{figure}

Because of the sign convention, the linking number turns out to be an
invariant of the link; it doesn't depend on the position of the link
in space.  See Figure~\ref{fig:linkinvar}, which shows invariance
under the one Reidemeister move (see Kauffman~\cite{kauffman:knotsphysics}) in
which the set of
$\A\B$-crossings changes.

Configuration space provides a nice way to understand
this combinatorial formula topologically.  Because of the type of
crossings that we count, it makes sense to consider the subset of
$\oC_2(\A\amalg \B)$ in which one point lies on $\A$ and the other on $\B$.
Since
$\A \cap \B = \emptyset$, this is isomorphic to $\oC_1(\A) \times
\oC_1(\B) \simeq \A \times \B$, which is a torus.

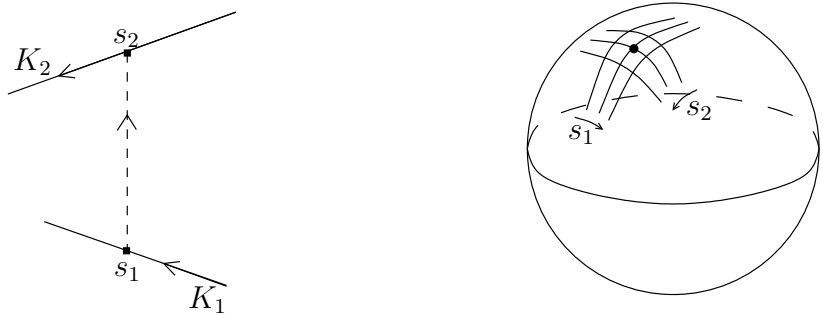
\begin{figure}
\centerline{\input{draws/crossing-sign.tex}}
\caption{The sign of a crossing is the sign of the local map from the
torus $A \times B$ to the sphere.  On the left are the sections of a
knot near this crossing.  Note that the view is not from
above, so it doesn't look like a crossing.  On the right are the
corresponding directions on the sphere.}
\label{fig:crossingsign}
\end{figure}

As in Section~\ref{sec:motivation}, we get a map $C^2(L): C^2(\A
\amalg \B) \rightarrow C_2(\RR^3)$. 
There is also a map\footnote{In fact,
a homotopy equivalence.}
$\phi_{12}: \oC_2(\RR^3) \rightarrow S^2$, given by
\begin{equation}
\phi_{12}(x_1, x_2) = \frac{x_2 - x_1}{|x_2 - x_1|}
\end{equation}
(the denominator does not vanish in $\oC_2(\RR^3)$).  Composing these
two, we get
a map $\phi_{12} \circ C^2(L)$ from the torus to the sphere.  The
``crossings'' we counted above were just those points on the torus
mapping to the vertical direction
 on the sphere.  If we orient both
the torus and the sphere, the sign of the
crossing defined above turns out to be the sign of the local mapping
from the torus to the sphere (i.e., whether it locally preserves or
reverses orientation), which will be non-degenerate for knots in
general position. See
Figure~\ref{fig:crossingsign}.

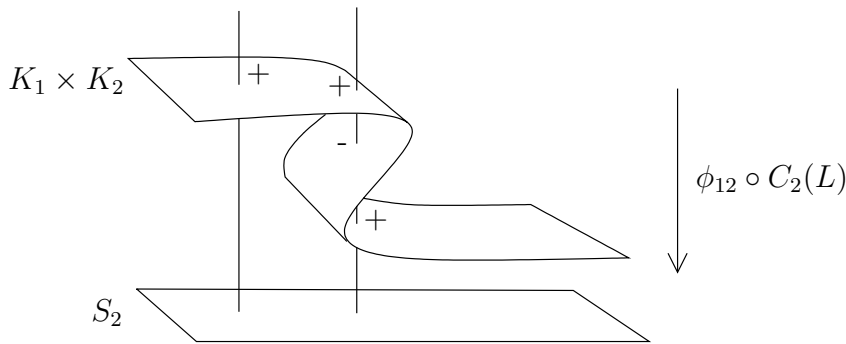
\begin{figure}
\centerline{\input{draws/fold.tex}}
\caption{The degree is independent of the point on the sphere used to
choose it.  As the point on the sphere moves from left to right, it
crosses a fold, introducing two new inverse images; but they're
counted with opposite signs.}
\label{fig:folding}
\end{figure}

In fact, this number is independent of the direction ($=$ point on the
sphere) that we choose, as long as it's a regular value.  See
Figure~\ref{fig:folding}.
This type of counting is well-known in topology; it is known as the
\emph{degree}, and can be defined for any map $f: M \rightarrow N$ between
oriented compact manifolds of the same dimension $n$.

To
define the degree rigourously, note that the highest cohomology of a
compact oriented
$n$-manifold $M$, $H^n(M)$, is always 1-dimensional and has a natural
generator $[\omega]$ (so that the integral of $\omega$ over the manifold
is 1); this is called the \emph{fundamental class} of $M$.  The map
$f: M \rightarrow N$ induces a map $f^*: H^n(N)
\rightarrow H^n(M)$,
taking the fundamental  class of $N$ to some multiple of the
fundamental class of $M$, which can be computed by integrating
$f^*\omega$ over $M$.
This multiple turns out to be exactly the
degree.  By Sard's theorem, we know that there is some regular point $x$
of $N$, i.e., points $x$ for which $f$ is a diffeomorphism in the
neighborhood of each
inverse image of $x$.
Since $M$ is compact, there are only a finite number of
inverse images of $x$.  If we take a representative $\omega$ of the
fundamental
class of $N$ that is sufficiently localized around $x$, then
the integral of $f^*\omega$ over $M$ is exactly the
degree as defined above.

But you need not pick a localized $\omega$.  In our case, of a map
from $S^1 \times S^1 \rightarrow S^2$,
it's natural to choose $\omega$ to be the
unique SO(3) invariant volume form on $S^2$,
\begin{align}
\omega 
       &= \frac{\epsilon_{\mu\nu\sigma}}{8\pi}\,
		\frac{x^\mu\,dx^\nu\,dx^\sigma}{|x|^3}
\end{align}
(using physicists' conventions: the components of
the 3-vector $x$ are $x^\mu$, $\mu = 1, 2, 3$ and repeated indices in
an expression are to be summed over, while $\epsilon_{\mu\nu\sigma}$
is the basic totally anti-symmetric tensor), then the integral of the
pullback of $\omega$ over the knot becomes
becomes
\begin{align}
\text{Linking number}
&= \int_{S^1_1} dx_1^\mu \,
   \int_{S^1_2} dx_2^\nu \,
   \frac{\epsilon_{\mu\nu\sigma}}{4\pi}\,\frac{(x_2 - x_1)^\sigma}
					{|x_2-x_1|^3} \\
 &= \int_{S^1 \times S^1} \theta_{12}
\end{align}
The first expression is again written in physicists' notation
(the $x_i^\mu$ are the components of the embedding of the $i$'th circle
into $\RR^3$).  In the second, I use notation from Bott and
Taubes~\cite{bt:knot}:
$\theta_{12} =  \circ C_2(K)^* \phi_{12}^* \omega$.

This is the simplest case of the integrals considered in this thesis.
(Although I won't consider links any further, the generalization of
these integrals to links is straightforward.)

\section{Back to Knots}

\begin{figure}
\centerline{\input{draws/cross-boundary.tex}}
\caption{Number of inverse images changes as you cross a boundary.}
\label{fig:crossboundary}
\end{figure}
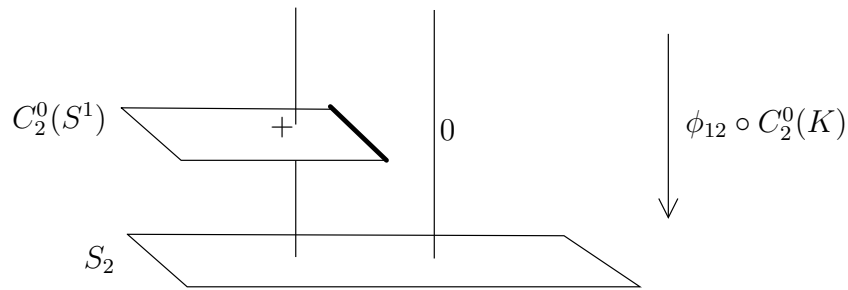

But this case is misleadingly simple.
Consider trying to modify
the linking invariant to get an invariant of a single knot.  As before, there
is a canonical map $\phi_{12} \circ \oC_2(K): \oC_2(S^1) \rightarrow
S^2$.  But now $\oC_2(S^1)$ is not
compact
(it is the torus $S^1
\times S^1$ minus the diagonal)
so there is no well-defined degree.
Alternatively, if you compactify $\oC_2(S^1)$ by adding a boundary, the
number of inverse images of a point in $S^2$
changes as you cross the image of the
boundary.  See Figure~\ref{fig:crossboundary}.

The resolution of this difficulty is to
take linear combinations of maps from spaces associated with the knot, each
with a boundary, in such a way that the boundaries exactly cancel out
in the target manifold.
The resulting ``degree'' will automatically be an invariant of knots.

In the case of $\oC_2(S^1)$, the necessary correction term comes from
a \emph{framing}
of the knot: a nowhere zero section of the normal bundle to the knot
(as imbedded in $\RR^3$).  (A framed knot can also be thought of as an
imbedded (oriented or two-sided) ribbon in $\RR^3$.)  With this
correction, you get 
the \emph{framing number} of the knot: the linking number of the knot
with a copy slightly displaced in the direction of the framing.

But it would
take us too far afield to describe this in more detail.
Instead, let's go to level 4, where there is an honest knot
invariant (rather than an invariant of framed knots).  That is,
consider $\oC_4(S^1)$.  Unlike $\oC_2(S^1)$, this is a disconnected space.
Pick the connected component in which the points 1, 2, 3, 4 appear in
that (cyclic) order around the circle.
There are now maps $\phi_{ij}: \oC_4(S^1) \rightarrow S^2$ (for $i \ne
j$, $1 \le i,j \le 4$), given by normalizing the vector $K(s_j) -
K(S_i)$.  To define a degree of any sort, we'll need a map from the
4-dimensional space $\oC_4(S^1)$ to another 4-dimensional manifold.  We
can find one such by considering, for instance, $\phi_{13} \times
\phi_{24}$.  (This is really the only possible choice using our
$\phi_{ij}$.  All of 1,2,3,4 must be represented, or the map is
degenerate, and if we have a term like $\phi_{12}$ (with two adjacent
points on the circle) we'll get a boundary just like that from
$\phi_{12}$ on $\oC_2(S^1)$, which  requires a framing-dependent
correction.)

The term to counter the boundary of the map $\phi_{13} \times
\phi_{24}$ turns out to be slightly more complicated.  It involves a
map from $\oC_{4;3}(\RR^3; K)$, the space of 4 points in $\RR^3$, with
the first 3 restricted to lie on the knot $K$.  From this
6-dimensional space we consider the map $\phi_{14} \times \phi_{24}
\times \phi_{34}: \oC_{4;3}(\RR^3; K) \rightarrow S^2 \times S^2
\times S^2$.  Then the resulting integral
\begin{equation}
\label{eq:v2}
\frac{1}{4}\int_{\oC_4} \theta_{13}\theta_{24} -
	\frac{1}{3}\int_{\oC_{4;3}} \theta_{14}\theta_{24}\theta_{34}
\end{equation}
converges and is a knot invariant.  (As before, $\theta_{ij} =
C_k(K)^* \circ
\phi_{ij}^*(\omega)$, with $\omega$ the normalized volume form on the
sphere.)

Raoul Bott and Cliff Taubes~\cite{bt:knot} proved this integral
is invariant.  (Earlier, Dror
Bar-Natan~\cite{bn:phd} had
proven invariance in another way.)
Their proof relies crucially on the  functorial
compactification $C_4(\RR^3)$ of the configuration space
$\oC_4(\RR^3)$, described in
detail in
Section~\ref{sec:configspace}.  This compactification yields a
manifold with corners, whose (codimension 1) faces correspond to
the points in a subset $B \subset \{1, 2, 3, 4\}$ of size $\ge 2$
approaching each other.  (Since the compactification is functorial,
the same description holds for the other spaces involved, $C_4(S^1)$
and $C_{4;3}(\RR^3;K)$.)

The proof is essentially Stokes' theorem: the variation of the
integral of a (closed)
form is the integral of the form on the boundary.  (This yields a
1-form on the space of knots $\knots$.)  See
Proposition~\ref{prop:dpushforward} for the precise statement.
Bott and Taubes showed that (a) either
integral restricted to any face in which more than two points come
together (as a 1-form on the space of knots) vanishes and (b) the
integrals restricted to the \emph{principal faces}, those in which
exactly two points come together, cancel each other out.

The principal faces are more interesting, so let's look at them first.
On the face of $C_4(S^1)$ in which points 1 and 2 come together, the
form $\theta_{13}\theta_{24}$ becomes
$\theta_{\{1,2\}3}\theta_{\{1,2\}4}$.  (The notation $\theta_{\{1,2\}3}$
indicates that $\theta_{13}$ is equal to $\theta_{23}$.  It also has
another meaning, relating to the coordinates on $C_4(S^1)$: on this
face, coordinates are the positions $s_{\{1,2\}}$, $s_3$, and $s_4$.
Thus this form can be considered a form on $C_{\{\{1,2\},3,4\}}(S^1)$.)
See Section~\ref{sec:configdescrip}.)  Similarly for the
other principal faces (with 2 and 3, 3 and 4, and 4 and 1 colliding)
of $C_4(S^1)$.  In $C_{4;3}(\RR^3; K)$, consider
the face in which points 1 and 4 come together.  On this face, the
form $\theta_{14}\theta_{24}\theta_{34}$ becomes
$\theta_{14}\theta_{2\{1,4\}}\theta_{3\{1,4\}}$.  But, on this face, the
map $\phi_{14}$ varies over $S^2$ without changing the directions to
any of the other points; so we can integrate over this $S^2$, yielding
the form $\theta_{2\{1,4\}}\theta_{3\{1,4\}}$ on
$C_{\{\{1,4\},2,3\}}(S^1)$.  Similarly for the faces with 2 and 4 and 3
and 4 colliding.  On the faces in which two points on the knot collide
(for instance, 1 and 2), the restriction of this form vanishes.

Thus each of the four principal faces of $C_4(S^1)$ contributes a form like
$\theta_{\{1,2\}3}\theta_{\{1,2\}4}$.  All of these yield the same
1-form on $\knots$.  Each of the three interesting principal faces of
$C_{4;3}(\RR^3; K)$ yields another form like this, again yielding the
same 1-form on $\knots$.  Because of the factors and signs in the
expression~\ref{eq:v2}, all these forms cancel.

This result was not news to the knot theory community.  This
compactification of configuration space had been known and used, for
instance by Axelrod and
Singer~\cite{Axelrod-Singer:ChernPert,Axelrod-Singer:ChernPertII}, who
proved existence and invariance of the similar integral invariants for
homology 3-spheres, and it was common knowledge that the techniques
worked for the knot invariants too.  But I know of no elementary
exposition.  The Bott-Taubes paper~\cite{bt:knot} is a first step;
unfortunately, it
omits description of the higher invariants, which can be treated with
the same methods.  This thesis attempts to rectify the situation.

The invariants are all constructed from configuration space integrals
of this basic type.  Namely, given a graph $\Gamma$ with a cycle $c$
representing the knot, the vertices of $\Gamma$ define a configuration
space $C_{V(\Gamma);V(c)}(\RR^3; K)$, and the edges of $\Gamma$ define
a map $\phi(\Gamma)$ to $(S^2)^{|e|}$ (after choosing an ordering of the
edges $e$ and a direction for each edge).  The corresponding integral
$I(\Gamma)$ is the integral over $C_{V(\Gamma);V(c)}(\RR^3; K)$ of the
pullback of the canonical form on $(S^2)^{|e|}$ (which is just the
product of the canonical forms on the $S^2$ factors).  Note that for
this to be well-defined, there must be an orientation on
$C_{V(\Gamma);V(c)}(\RR^3; K)$.  An ordering of the vertices of
$\Gamma$ gives  such an orientation.



Chapter~\ref{sec:vassiliev} defines the Vassiliev invariants and
introduces weight diagrams, with a plausibility argument why they
should yield invariant integrals.  Chapter~\ref{sec:configspace}
introduces the functorial compactification of configuration space that
is essential to removing the analytic difficulties involved.
Chapter~\ref{sec:tying} ties up loose ends, completing the description
of the invariants, showing that the principle faces cancel, and
describing the tinkertoy diagram formula.
Chapter~\ref{sec:nastyface} shows that the contributions to the
variation of the integral from most
hidden or non-principle faces vanish.  The remaining faces (the
anomalous faces) are the subject of Chapter~\ref{sec:anomaly}.

There are two appendices, describing in more detail the signs involved
in relating weight diagrams and configuration space integrals.
Appendix~\ref{sec:linalg} gives some linear algebra preliminaries for
the main statements in Appendix~\ref{sec:graph}.

\section{Acknowledgements}

I would like to thank 
Ezra Getzler,
Michael Hutchings,
Jose Labastida,
Justin Roberts,
James Stasheff,
Clifford Taubes, and
Zhenghan Wang
for helpful conversations.  Dror Bar-Natan has helped straighten out
my confusions in many ways.  And Raoul Bott, my advisor,
has been tremendously helpful throughout countless hours of
discussion.

%% file: draws/link-sign.tex
\begin{picture}(0,0)%
\epsfig{file=draws/link-sign.pstex}%
\end{picture}%
\setlength{\unitlength}{3947sp}%
\begingroup\makeatletter\ifx\SetFigFont\undefined
\def\x#1#2#3#4#5#6#7\relax{\def\x{#1#2#3#4#5#6}}%
\expandafter\x\fmtname xxxxxx\relax \def\y{splain}%
\ifx\x\y   
\gdef\SetFigFont#1#2#3{%
  \ifnum #1<17\tiny\else \ifnum #1<20\small\else
  \ifnum #1<24\normalsize\else \ifnum #1<29\large\else
  \ifnum #1<34\Large\else \ifnum #1<41\LARGE\else
     \huge\fi\fi\fi\fi\fi\fi
  \csname #3\endcsname}%
\else
\gdef\SetFigFont#1#2#3{\begingroup
  \count@#1\relax \ifnum 25<\count@\count@25\fi
  \def\x{\endgroup\@setsize\SetFigFont{#2pt}}%
  \expandafter\x
    \csname \romannumeral\the\count@ pt\expandafter\endcsname
    \csname @\romannumeral\the\count@ pt\endcsname
  \csname #3\endcsname}%
\fi
\fi\endgroup
\begin{picture}(4184,1528)(1406,-2548)
\put(4389,-1215){\makebox(0,0)[lb]{\smash{\SetFigFont{12}{14.4}{rm}$K_2$}}}
\put(5512,-1281){\makebox(0,0)[lb]{\smash{\SetFigFont{12}{14.4}{rm}$K_1$}}}
\put(1715,-1215){\makebox(0,0)[lb]{\smash{\SetFigFont{12}{14.4}{rm}$K_2$}}}
\put(2837,-1281){\makebox(0,0)[lb]{\smash{\SetFigFont{12}{14.4}{rm}$K_1$}}}
\put(2416,-1831){\makebox(0,0)[lb]{\smash{\SetFigFont{12}{14.4}{rm}$-1$}}}
\put(5071,-1861){\makebox(0,0)[lb]{\smash{\SetFigFont{12}{14.4}{rm}$+1$}}}
\end{picture}

%% file: draws/link-Hopf.tex
\begin{picture}(0,0)%
\epsfig{file=draws/link-Hopf.pstex}%
\end{picture}%
\setlength{\unitlength}{3947sp}%
\begingroup\makeatletter\ifx\SetFigFont\undefined
\def\x#1#2#3#4#5#6#7\relax{\def\x{#1#2#3#4#5#6}}%
\expandafter\x\fmtname xxxxxx\relax \def\y{splain}%
\ifx\x\y   
\gdef\SetFigFont#1#2#3{%
  \ifnum #1<17\tiny\else \ifnum #1<20\small\else
  \ifnum #1<24\normalsize\else \ifnum #1<29\large\else
  \ifnum #1<34\Large\else \ifnum #1<41\LARGE\else
     \huge\fi\fi\fi\fi\fi\fi
  \csname #3\endcsname}%
\else
\gdef\SetFigFont#1#2#3{\begingroup
  \count@#1\relax \ifnum 25<\count@\count@25\fi
  \def\x{\endgroup\@setsize\SetFigFont{#2pt}}%
  \expandafter\x
    \csname \romannumeral\the\count@ pt\expandafter\endcsname
    \csname @\romannumeral\the\count@ pt\endcsname
  \csname #3\endcsname}%
\fi
\fi\endgroup
\begin{picture}(3937,1639)(1576,-2688)
\put(4560,-1313){\makebox(0,0)[lb]{\smash{\SetFigFont{12}{14.4}{rm}$K_2$}}}
\put(5513,-1792){\makebox(0,0)[lb]{\smash{\SetFigFont{12}{14.4}{rm}Linking $\# = +1$}}}
\put(1576,-1336){\makebox(0,0)[lb]{\smash{\SetFigFont{12}{14.4}{rm}$K_1$}}}
\put(3406,-2341){\makebox(0,0)[b]{\smash{\SetFigFont{12}{14.4}{rm}+}}}
\end{picture}

%% file: draws/link-inv.tex
\begin{picture}(0,0)%
\epsfig{file=draws/link-inv.pstex}%
\end{picture}%
\setlength{\unitlength}{3947sp}%
\begingroup\makeatletter\ifx\SetFigFont\undefined
\def\x#1#2#3#4#5#6#7\relax{\def\x{#1#2#3#4#5#6}}%
\expandafter\x\fmtname xxxxxx\relax \def\y{splain}%
\ifx\x\y   
\gdef\SetFigFont#1#2#3{%
  \ifnum #1<17\tiny\else \ifnum #1<20\small\else
  \ifnum #1<24\normalsize\else \ifnum #1<29\large\else
  \ifnum #1<34\Large\else \ifnum #1<41\LARGE\else
     \huge\fi\fi\fi\fi\fi\fi
  \csname #3\endcsname}%
\else
\gdef\SetFigFont#1#2#3{\begingroup
  \count@#1\relax \ifnum 25<\count@\count@25\fi
  \def\x{\endgroup\@setsize\SetFigFont{#2pt}}%
  \expandafter\x
    \csname \romannumeral\the\count@ pt\expandafter\endcsname
    \csname @\romannumeral\the\count@ pt\endcsname
  \csname #3\endcsname}%
\fi
\fi\endgroup
\begin{picture}(2844,1373)(1193,-1378)
\put(3523,-200){\makebox(0,0)[lb]{\smash{\SetFigFont{12}{14.4}{rm}$K_1$}}}
\put(4001,-200){\makebox(0,0)[lb]{\smash{\SetFigFont{12}{14.4}{rm}$K_2$}}}
\put(1236,-200){\makebox(0,0)[lb]{\smash{\SetFigFont{12}{14.4}{rm}$K_1$}}}
\put(1903,-200){\makebox(0,0)[lb]{\smash{\SetFigFont{12}{14.4}{rm}$K_2$}}}
\put(2780,-885){\makebox(0,0)[b]{\smash{\SetFigFont{12}{14.4}{rm}$\Longrightarrow$}}}
\put(2780,-610){\makebox(0,0)[b]{\smash{\SetFigFont{12}{14.4}{rm}Reidemeister move}}}
\put(2780,-404){\makebox(0,0)[b]{\smash{\SetFigFont{12}{14.4}{rm}Type II}}}
\put(3946,-556){\makebox(0,0)[b]{\smash{\SetFigFont{12}{14.4}{rm}+}}}
\put(3961,-1171){\makebox(0,0)[b]{\smash{\SetFigFont{12}{14.4}{rm}-}}}
\end{picture}

%% file: draws/crossing-sign.tex
\begin{picture}(0,0)%
\epsfig{file=draws/crossing-sign.pstex}%
\end{picture}%
\setlength{\unitlength}{3947sp}%
\begingroup\makeatletter\ifx\SetFigFont\undefined
\def\x#1#2#3#4#5#6#7\relax{\def\x{#1#2#3#4#5#6}}%
\expandafter\x\fmtname xxxxxx\relax \def\y{splain}%
\ifx\x\y   
\gdef\SetFigFont#1#2#3{%
  \ifnum #1<17\tiny\else \ifnum #1<20\small\else
  \ifnum #1<24\normalsize\else \ifnum #1<29\large\else
  \ifnum #1<34\Large\else \ifnum #1<41\LARGE\else
     \huge\fi\fi\fi\fi\fi\fi
  \csname #3\endcsname}%
\else
\gdef\SetFigFont#1#2#3{\begingroup
  \count@#1\relax \ifnum 25<\count@\count@25\fi
  \def\x{\endgroup\@setsize\SetFigFont{#2pt}}%
  \expandafter\x
    \csname \romannumeral\the\count@ pt\expandafter\endcsname
    \csname @\romannumeral\the\count@ pt\endcsname
  \csname #3\endcsname}%
\fi
\fi\endgroup
\begin{picture}(5199,1908)(1025,-1962)
\put(4727,-903){\makebox(0,0)[b]{\smash{\SetFigFont{12}{14.4}{rm}$s_1$}}}
\put(5468,-748){\makebox(0,0)[b]{\smash{\SetFigFont{12}{14.4}{rm}$s_2$}}}
\put(1876,-1762){\makebox(0,0)[b]{\smash{\SetFigFont{12}{14.4}{rm}$s_1$}}}
\put(1868,-295){\makebox(0,0)[b]{\smash{\SetFigFont{12}{14.4}{rm}$s_2$}}}
\put(2380,-1962){\makebox(0,0)[b]{\smash{\SetFigFont{12}{14.4}{rm}$K_1$}}}
\put(1279,-467){\makebox(0,0)[b]{\smash{\SetFigFont{12}{14.4}{rm}$K_2$}}}
\end{picture}

%% file: draws/fold.tex
\begin{picture}(0,0)%
\epsfig{file=draws/fold.pstex}%
\end{picture}%
\setlength{\unitlength}{3947sp}%
\begingroup\makeatletter\ifx\SetFigFont\undefined
\def\x#1#2#3#4#5#6#7\relax{\def\x{#1#2#3#4#5#6}}%
\expandafter\x\fmtname xxxxxx\relax \def\y{splain}%
\ifx\x\y   
\gdef\SetFigFont#1#2#3{%
  \ifnum #1<17\tiny\else \ifnum #1<20\small\else
  \ifnum #1<24\normalsize\else \ifnum #1<29\large\else
  \ifnum #1<34\Large\else \ifnum #1<41\LARGE\else
     \huge\fi\fi\fi\fi\fi\fi
  \csname #3\endcsname}%
\else
\gdef\SetFigFont#1#2#3{\begingroup
  \count@#1\relax \ifnum 25<\count@\count@25\fi
  \def\x{\endgroup\@setsize\SetFigFont{#2pt}}%
  \expandafter\x
    \csname \romannumeral\the\count@ pt\expandafter\endcsname
    \csname @\romannumeral\the\count@ pt\endcsname
  \csname #3\endcsname}%
\fi
\fi\endgroup
\begin{picture}(4320,2124)(2956,-2473)
\put(4531,-826){\makebox(0,0)[b]{\smash{\SetFigFont{12}{14.4}{rm}+}}}
\put(5266,-1741){\makebox(0,0)[b]{\smash{\SetFigFont{12}{14.4}{rm}+}}}
\put(5041,-901){\makebox(0,0)[b]{\smash{\SetFigFont{12}{14.4}{rm}+}}}
\put(5056,-1261){\makebox(0,0)[b]{\smash{\SetFigFont{12}{14.4}{rm}-}}}
\put(7276,-1486){\makebox(0,0)[lb]{\smash{\SetFigFont{12}{14.4}{rm}$\phi_{12} \circ C_2(L)$}}}
\put(2956,-871){\makebox(0,0)[lb]{\smash{\SetFigFont{12}{14.4}{rm}$K_1 \times K_2$}}}
\put(3481,-2326){\makebox(0,0)[lb]{\smash{\SetFigFont{12}{14.4}{rm}$S_2$}}}
\end{picture}

%% file: draws/cross-boundary.tex
\begin{picture}(0,0)%
\epsfig{file=draws/cross-boundary.pstex}%
\end{picture}%
\setlength{\unitlength}{3947sp}%
\begingroup\makeatletter\ifx\SetFigFont\undefined
\def\x#1#2#3#4#5#6#7\relax{\def\x{#1#2#3#4#5#6}}%
\expandafter\x\fmtname xxxxxx\relax \def\y{splain}%
\ifx\x\y   
\gdef\SetFigFont#1#2#3{%
  \ifnum #1<17\tiny\else \ifnum #1<20\small\else
  \ifnum #1<24\normalsize\else \ifnum #1<29\large\else
  \ifnum #1<34\Large\else \ifnum #1<41\LARGE\else
     \huge\fi\fi\fi\fi\fi\fi
  \csname #3\endcsname}%
\else
\gdef\SetFigFont#1#2#3{\begingroup
  \count@#1\relax \ifnum 25<\count@\count@25\fi
  \def\x{\endgroup\@setsize\SetFigFont{#2pt}}%
  \expandafter\x
    \csname \romannumeral\the\count@ pt\expandafter\endcsname
    \csname @\romannumeral\the\count@ pt\endcsname
  \csname #3\endcsname}%
\fi
\fi\endgroup
\begin{picture}(4230,1779)(3046,-2473)
\put(4666,-1501){\makebox(0,0)[lb]{\smash{\SetFigFont{12}{14.4}{rm}+}}}
\put(5731,-1531){\makebox(0,0)[lb]{\smash{\SetFigFont{12}{14.4}{rm}0}}}
\put(7276,-1486){\makebox(0,0)[lb]{\smash{\SetFigFont{12}{14.4}{rm}$\phi_{12} \circ C_2^0(K)$}}}
\put(3496,-2341){\makebox(0,0)[lb]{\smash{\SetFigFont{12}{14.4}{rm}$S_2$}}}
\put(3046,-1456){\makebox(0,0)[lb]{\smash{\SetFigFont{12}{14.4}{rm}$C_2^0(S^1)$}}}
\end{picture}

%% file: vassiliev.tex
\chapter{Vassiliev Invariants and Weight Diagrams}\label{sec:vassiliev}

This chapter is a brief introduction to Vassiliev
invariants and the combinatorics of weight diagrams,
which are basic to the construction of
configuration space integrals.  The reader is referred to Bar-Natan's
excellent 
paper~\cite{bn:vassiliev} for a more in-depth treatment.

Any knot invariant $V: \knots
\rightarrow \RR$ can
be extended to an invariant $V^{(k)}$ of immersed
circles with exactly $k$ transversal self-intersections inductively,
using
\begin{align}
V^{(0)} &= V\\
\label{eq:crossing}
V^{(k)}(\mathcenter{\input{draws/intersect.tex}}) &= V^{(k-1)}(\mathcenter{\input{draws/overcross.tex}}) - V^{(k-1)}(\mathcenter{\input{draws/undercross.tex}})
\end{align}
Equation~\ref{eq:crossing} is to be interpreted a local statement.  The
loose ends are connected by some paths, possibly with crossings.  Note
that this definition depends on the orientation on $\RR^3$, but
not the plane projection or the orientation on $S^1$.  This is a good
definition, since (expanding out) it expresses $V^{(k)}$ on a ``knot''
with $k$ transversal self-intersections as a signed sum of the values
of $V$ on the $2^k$
ways of resolving the singularity, which won't depend on the order of
expansion.

\begin{definition}[The Birman-Lin Condition]
A \emph{Vassiliev invariant}~\cite{Vas:CohKnot,Vas:Book,bn:vassiliev}
of degree $k$
is a knot invariant that vanishes on any knot with $k+1$ double
points, when extended as above.
\end{definition}

One way to think of these Vassiliev invariants is as invariants of
polynomial type.
Taking the difference of the value of the knot invariant between the two
sides of a crossing can be seen as ``differentiation'' of a sort.
A polynomial is a function that becomes 0 if you differentiate it
enough, similar to the definition of a Vassiliev invariant above.

Directly from the definition, we can see that the invariants we'll
construct will be Vassiliev invariants.

\begin{proposition}\label{prop:integralvas}
A sum of 
canonical configuration space integrals of degree $k$ (involving a
configuration space with
$2k$ points) that is a knot invariant is a
Vassiliev invariant of degree $k$.
\end{proposition}

The proof is in Section~\ref{sec:integralvas}.  The essential idea is
that there can only be a difference between an undercrossing and an
overcrossing if there is at
least one point on each strand very near the crossing and a propogator
connecting them.

\begin{figure}
\centerline{\input{draws/chord-diag.tex}}
\caption{A chord diagram corresponding to a knot with 2 double
points.}
\label{fig:chorddiag}
\end{figure}
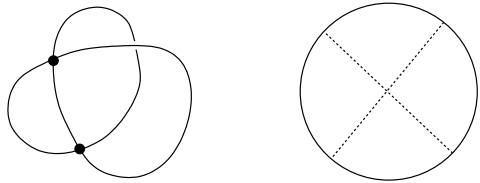

If you differentiate a Vassiliev invariant of degree $k$ just $k$
times, you get the analog of a constant.  This is an invariant of
knots with $k$ double points
that doesn't change as the knot passes through itself: this is exactly
the condition that it vanishes on knots with $k+1$ double points.
Thus the invariant evaluated on a knot with $k$ double points doesn't
depend on the embedding in space, but only on the topological
structure of the image.
This is completely determined by 
labelling the points of self-intersection of the knot and recording
the labels of the points as you travel around the circle in the
increasing coordinate.  You can get a graphical representation (called
a \emph{chord diagram}) by
drawing the circle as a circle and connecting points on the circle
that get mapped to the same point in $\RR^3$, as in
Figure~\ref{fig:chorddiag}.

Each Vassiliev invariant of degree $k$ therefore determines a function
on chord diagrams with $k$ chords.  We can evaluate this function
explicitly for configuration space integrals.  By the proof of
Proposition~\ref{prop:integralvas}, the only configuration space
integrals of degree $k$ that don't vanish when evaluated on a knot
with $k$ double points are those that involve $2k$ points, all on the
knot, connected by propogators in pairs that are near the double
points.  In particular, the graph of the propogators must be
isomorphic the the graph of the chord diagram.
Since the integral of a basic propogator near a double point is
$\pm1$, the value of the configuration space integral on this knot
with double points will be (up to sign) the number of such
isomorphisms.  See Section~\ref{sec:integralvas} for more details.

\begin{figure}
\centerline{\input{draws/4-term.tex}}
\caption{The 4 term relation, imbedded in $\RR^3$ (left) and
unfolded into a chord diagram (right).  The dotted sections of the
circle are portions that may have other chords attached (possibly
crossing those drawn).}
\label{fig:4T}
\end{figure}
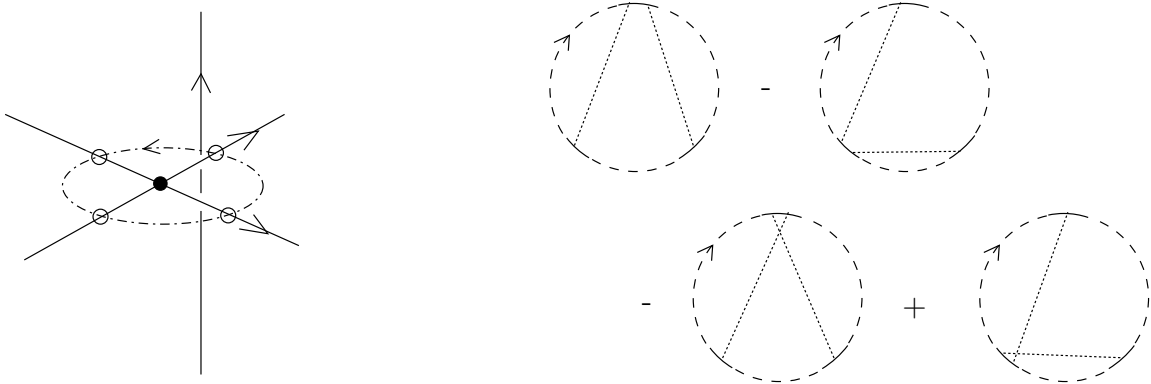

Not all functions on chord diagrams can
be achieved as derivatives of knot invariants.  In particular, imagine
evaluating the invariant on a knot
with $k-1$ double points and moving one strand in a small circle
around one of the other double points, as on the left in
Figure~\ref{fig:4T}.
The difference as you pass through each strand is the value of the
invariant on a certain chord diagram.  Since you end up at the same
knot you started at, the sum of the value on the four diagrams
produced must be 0.  Unfolding the relation, you get the sum of
diagrams on the right in Figure~\ref{fig:4T}, call the \emph{four
term} or $4T$ relation.  Functions on chord diagrams satisfying the
$4T$ relation are called \emph{weight diagrams}.
One of the consequences of this thesis is
that any weight diagram
can be ``integrated'' to yield an invariant of knots.  This was proved
earlier by Kontsevich using a different type of integral.
See~\cite{bn:vassiliev} for a description.

The number of weight diagrams of a given degree is an important open
question.
Bar-Natan~\cite{bn:vassiliev} has
done extensive computer computations, computing all weight diagrams up
to degree 9.  There seem to be a fair number of them, probably growing
exponentially, although the difficulty of computation grows as a
considerably
faster exponential.

It's an amazing 
 fact that any weight
diagram can be extended to a function on
diagrams with internal vertices of a certain type.  Namely, a
\emph{trivalent graph} (termed a Chinese Character Diagram
by Bar-Natan~\cite{bn:vassiliev}) is a connected graph $\Gamma$ with a
directed cycle $c$
(with distinct vertices) and all vertices trivalent.  Vertices that do
not meet $c$ are the \emph{internal} (or free or interaction)
vertices.  The value of a weight
diagram on a trivalent graph will depend on an orientation
on the edges incident to each internal vertex.  Since the vertices are
trivalent, this is a cylic order.

\begin{proposition}[Bar-Natan~\cite{bn:vassiliev},Theorem~6]
\label{prop:weightdiag}
A weight diagram $w$ extends uniquely to a function on trivalent graphs
with an orientation on the edges incident to each vertex
satisfying the
$STU$ relation
\begin{equation}
\mathcenter{\input{draws/STU.tex}}
\end{equation}
\end{proposition}

\begin{figure}
\centerline{\input{draws/sillycase.tex}}
\caption{The remaining case of Proposition~\ref{prop:weightdiag}.
There is an internal blob $B$ connected to a trivalent vertex and
thence to a circle, but no other trivalent vertices.}
\label{fig:sillycase}
\end{figure}
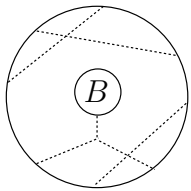

\begin{proof}
It's clear that (since trivalent graphs are connected) the $STU$
relation defines inductively the value of $w$ on any
trivalent graph, as long as the
 definition does not introduce
inconsistency.

The $4T$ relation is a difference of two
$STU$ relations involving graphs with zero or one internal vertices.
Because of this, the value of $w$ on a graph with exactly
one internal vertex is well-defined.  The proof then proceeds by induction.
Suppose $w$ has been defined on graphs
with $k\ge1$ or fewer internal vertices.  Consider $w(\Gamma_{k+1})$ for a
graph $\Gamma_{k+1}$ with $k+1$ vertices.  This can be defined
from the value of $w$ on graphs with $k$ internal vertices
via $STU$ on any of
the edges that connects an internal vertex of $\Gamma_{k+1}$ to the circle.
If two of these edges connect to different internal vertices, 
then the corresponding two potential expressions for $w(\Gamma_{k+1})$
can be expanded
into a sum of $w(\Gamma'_{k-1})$ for graphs $\Gamma'_{k-1}$ with $k-1$
vertices via
$STU$ at the next lower level on the other edge.  The resulting graphs
$\Gamma'$ are the same, so the two potential expressions for
$w(\Gamma_{k+1})$ are consistent.  If two edges of $\Gamma_{k+1}$
connect the
same internal vertex to the circle but there is a third edge
connecting a different internal vertex to the circle, we can apply
this procedure in two steps to deduce consistency.  There is one
remaing case, as in Figure~\ref{fig:sillycase}.  Since the value of
the weight diagram on trivalent diagrams like this will turn out to be
irrelevant, I defer to Bar-Natan~\cite{bn:vassiliev} for the completion of the
proof.
\end{proof}

Internal vertices arise naturally in the Feynmann diagrams for
Chern-Simons perturbation theory, where they represent the
interactions of the theory.  But they are still a mystery from other
points of view.  $w(\Gamma)$ has not yet been given a sensible
interpretation in terms of knots when $\Gamma$ has internal vertices.
Bott and Taubes~\cite{bt:knot} have suggested that, from the
differential point of view they are related to
Sullivan's minimal model for the rational homotopy of $\oC_n(\RR^3)$,
but this has not yet been explained clearly.  In any case, we have the
following theorem. $I(\Gamma)$ is, as in the introduction, the
configuration space
integral corresponding to the trivalent graph $\Gamma$ (with a sign yet to be
defined).

\begin{theorem}\label{thm:mainthm}
For any weight diagram $w$, the following sum of integrals is a knot
invariant:
\begin{equation}
\label{eq:mainthm}
I(w) = \sum_{\substack{\mathrm{isomorphism\ classes}\\
	\mathrm{of\ graphs\ }\Gamma}}
\frac{w(\Gamma) I(\Gamma)}{|\mathrm{Aut}(\Gamma)|}
+ \mathrm{(anomaly\ correction\ term)}
\end{equation}
\end{theorem}

The expression~\ref{eq:mainthm}, without the correction term, is
exactly the integral that would be
written down in perturbation theory, but it has been decoupled from
all physics considerations here.

The proof is spread out over the remainder of the paper.
The basic idea is that the differential form extends smoothly to the
compactification of configuration space.  In this context, the
variation of the
integral as the knot varies can be written (up to sign) as the sum of the
restrictions of the integral to the faces of configuration space.
(See Proposition~\ref{prop:dpushforward}.)

The most interesting faces are the principle faces, which correspond
to exactly two points (connected by either a segment of the knot or a
propogator) approaching each other.  The $STU$ relation expresses the
cancellation of faces in which two vertices on the knot or a vertex on
the knot and an internal vertex approach each other.  The remaining
principle faces (two internal vertices colliding) cancel by the
following proposition:

\begin{proposition}[Bar-Natan~\cite{bn:vassiliev}, Theorem~6]
Every weight diagram $w$ extended to graphs satisfies the
$IHX$ relation:
\begin{equation}
\mathcenter{\input{draws/IHX.tex}}
\end{equation}
\end{proposition}

I will not repeat the proof here.  It reduces to the case when one of
the edges is connected to the circle, and in this case you get a sum
of terms very much like the Jacobi identity for matrices.

The factor $|\mathrm{Aut}(\Gamma)|$, the size of the automorphism
group of the graph $\Gamma$, arises for the same reason a similar
factor appeared for the evaluation of an integral on a knot with
double points above.  A different explanation is that the number of
edges like the any given edge in $\Gamma$, and hence the number of
faces of configuration space
whose variation is the contraction of that edge, is usually the number
of automorphisms of $\Gamma$.

The remaining faces of configuration space correspond to some number
$k \ge 3$ of vertices approaching each other.  Most of these are shown
to vanish in Chapter~\ref{sec:nastyface}.  The remaining faces are the
\emph{anomalous faces}: those faces in which the $k$ degenerating
vertices are not connected to any other vertex by a propogator.  (This
includes the case in which all the vertices degenerate.)  The
variations resulting from these faces is described at length in
Chapter~\ref{sec:anomaly}, and is the source of the
(potential) 
correction term in equation~\ref{eq:mainthm}.

%% file: draws/chord-diag.tex
\begin{picture}(0,0)%
\epsfig{file=draws/chord-diag.pstex}%
\end{picture}%
\setlength{\unitlength}{3947sp}%
\begingroup\makeatletter\ifx\SetFigFont\undefined
\def\x#1#2#3#4#5#6#7\relax{\def\x{#1#2#3#4#5#6}}%
\expandafter\x\fmtname xxxxxx\relax \def\y{splain}%
\ifx\x\y   
\gdef\SetFigFont#1#2#3{%
  \ifnum #1<17\tiny\else \ifnum #1<20\small\else
  \ifnum #1<24\normalsize\else \ifnum #1<29\large\else
  \ifnum #1<34\Large\else \ifnum #1<41\LARGE\else
     \huge\fi\fi\fi\fi\fi\fi
  \csname #3\endcsname}%
\else
\gdef\SetFigFont#1#2#3{\begingroup
  \count@#1\relax \ifnum 25<\count@\count@25\fi
  \def\x{\endgroup\@setsize\SetFigFont{#2pt}}%
  \expandafter\x
    \csname \romannumeral\the\count@ pt\expandafter\endcsname
    \csname @\romannumeral\the\count@ pt\endcsname
  \csname #3\endcsname}%
\fi
\fi\endgroup
\begin{picture}(2968,1126)(3432,-2814)
\put(4981,-2371){\makebox(0,0)[b]{\smash{\SetFigFont{17}{20.4}{rm}~}}}
\end{picture}

%% file: draws/4-term.tex
\begin{picture}(0,0)%
\epsfig{file=draws/4-term.pstex}%
\end{picture}%
\setlength{\unitlength}{3947sp}%
\begingroup\makeatletter\ifx\SetFigFont\undefined
\def\x#1#2#3#4#5#6#7\relax{\def\x{#1#2#3#4#5#6}}%
\expandafter\x\fmtname xxxxxx\relax \def\y{splain}%
\ifx\x\y   
\gdef\SetFigFont#1#2#3{%
  \ifnum #1<17\tiny\else \ifnum #1<20\small\else
  \ifnum #1<24\normalsize\else \ifnum #1<29\large\else
  \ifnum #1<34\Large\else \ifnum #1<41\LARGE\else
     \huge\fi\fi\fi\fi\fi\fi
  \csname #3\endcsname}%
\else
\gdef\SetFigFont#1#2#3{\begingroup
  \count@#1\relax \ifnum 25<\count@\count@25\fi
  \def\x{\endgroup\@setsize\SetFigFont{#2pt}}%
  \expandafter\x
    \csname \romannumeral\the\count@ pt\expandafter\endcsname
    \csname @\romannumeral\the\count@ pt\endcsname
  \csname #3\endcsname}%
\fi
\fi\endgroup
\begin{picture}(7202,2405)(1789,-3954)
\put(6586,-2146){\makebox(0,0)[b]{\smash{\SetFigFont{12}{14.4}{rm}-}}}
\put(5836,-3496){\makebox(0,0)[b]{\smash{\SetFigFont{12}{14.4}{rm}-}}}
\put(7516,-3526){\makebox(0,0)[b]{\smash{\SetFigFont{12}{14.4}{rm}+}}}
\end{picture}

%% file: draws/STU.tex
\begin{picture}(0,0)%
\epsfig{file=draws/STU.pstex}%
\end{picture}%
\setlength{\unitlength}{3947sp}%
\begingroup\makeatletter\ifx\SetFigFont\undefined
\def\x#1#2#3#4#5#6#7\relax{\def\x{#1#2#3#4#5#6}}%
\expandafter\x\fmtname xxxxxx\relax \def\y{splain}%
\ifx\x\y   
\gdef\SetFigFont#1#2#3{%
  \ifnum #1<17\tiny\else \ifnum #1<20\small\else
  \ifnum #1<24\normalsize\else \ifnum #1<29\large\else
  \ifnum #1<34\Large\else \ifnum #1<41\LARGE\else
     \huge\fi\fi\fi\fi\fi\fi
  \csname #3\endcsname}%
\else
\gdef\SetFigFont#1#2#3{\begingroup
  \count@#1\relax \ifnum 25<\count@\count@25\fi
  \def\x{\endgroup\@setsize\SetFigFont{#2pt}}%
  \expandafter\x
    \csname \romannumeral\the\count@ pt\expandafter\endcsname
    \csname @\romannumeral\the\count@ pt\endcsname
  \csname #3\endcsname}%
\fi
\fi\endgroup
\begin{picture}(4801,1092)(2993,-3886)
\put(3601,-3886){\makebox(0,0)[b]{\smash{\SetFigFont{12}{14.4}{rm}S}}}
\put(5401,-3886){\makebox(0,0)[b]{\smash{\SetFigFont{12}{14.4}{rm}T}}}
\put(7201,-3871){\makebox(0,0)[b]{\smash{\SetFigFont{12}{14.4}{rm}U}}}
\put(4441,-3211){\makebox(0,0)[b]{\smash{\SetFigFont{12}{14.4}{rm}=}}}
\put(6286,-3226){\makebox(0,0)[b]{\smash{\SetFigFont{12}{14.4}{rm}-}}}
\end{picture}

%% file: draws/sillycase.tex
\begin{picture}(0,0)%
\epsfig{file=draws/sillycase.pstex}%
\end{picture}%
\setlength{\unitlength}{3947sp}%
\begingroup\makeatletter\ifx\SetFigFont\undefined
\def\x#1#2#3#4#5#6#7\relax{\def\x{#1#2#3#4#5#6}}%
\expandafter\x\fmtname xxxxxx\relax \def\y{splain}%
\ifx\x\y   
\gdef\SetFigFont#1#2#3{%
  \ifnum #1<17\tiny\else \ifnum #1<20\small\else
  \ifnum #1<24\normalsize\else \ifnum #1<29\large\else
  \ifnum #1<34\Large\else \ifnum #1<41\LARGE\else
     \huge\fi\fi\fi\fi\fi\fi
  \csname #3\endcsname}%
\else
\gdef\SetFigFont#1#2#3{\begingroup
  \count@#1\relax \ifnum 25<\count@\count@25\fi
  \def\x{\endgroup\@setsize\SetFigFont{#2pt}}%
  \expandafter\x
    \csname \romannumeral\the\count@ pt\expandafter\endcsname
    \csname @\romannumeral\the\count@ pt\endcsname
  \csname #3\endcsname}%
\fi
\fi\endgroup
\begin{picture}(1164,1159)(2164,-2498)
\put(2756,-1942){\makebox(0,0)[b]{\smash{\SetFigFont{12}{14.4}{rm}$B$}}}
\end{picture}

%% file: draws/IHX.tex
\begin{picture}(0,0)%
\epsfig{file=draws/IHX.pstex}%
\end{picture}%
\setlength{\unitlength}{3947sp}%
\begingroup\makeatletter\ifx\SetFigFont\undefined
\def\x#1#2#3#4#5#6#7\relax{\def\x{#1#2#3#4#5#6}}%
\expandafter\x\fmtname xxxxxx\relax \def\y{splain}%
\ifx\x\y   
\gdef\SetFigFont#1#2#3{%
  \ifnum #1<17\tiny\else \ifnum #1<20\small\else
  \ifnum #1<24\normalsize\else \ifnum #1<29\large\else
  \ifnum #1<34\Large\else \ifnum #1<41\LARGE\else
     \huge\fi\fi\fi\fi\fi\fi
  \csname #3\endcsname}%
\else
\gdef\SetFigFont#1#2#3{\begingroup
  \count@#1\relax \ifnum 25<\count@\count@25\fi
  \def\x{\endgroup\@setsize\SetFigFont{#2pt}}%
  \expandafter\x
    \csname \romannumeral\the\count@ pt\expandafter\endcsname
    \csname @\romannumeral\the\count@ pt\endcsname
  \csname #3\endcsname}%
\fi
\fi\endgroup
\begin{picture}(6024,1766)(2389,-3211)
\put(4201,-2236){\makebox(0,0)[b]{\smash{\SetFigFont{12}{14.4}{rm}=}}}
\put(6601,-2236){\makebox(0,0)[b]{\smash{\SetFigFont{12}{14.4}{rm}-}}}
\put(3001,-3211){\makebox(0,0)[b]{\smash{\SetFigFont{12}{14.4}{rm}$I$}}}
\put(4201,-3211){\makebox(0,0)[b]{\smash{\SetFigFont{12}{14.4}{rm}=}}}
\put(5401,-3211){\makebox(0,0)[b]{\smash{\SetFigFont{12}{14.4}{rm}$H$}}}
\put(6601,-3211){\makebox(0,0)[b]{\smash{\SetFigFont{12}{14.4}{rm}-}}}
\put(7801,-3211){\makebox(0,0)[b]{\smash{\SetFigFont{12}{14.4}{rm}$X$}}}
\end{picture}

%% file: configspace.tex
\chapter{Configuration space and its Compactification}
\label{sec:configspace}

\begin{definition}
The \emph{configuration space} $\oC_A(M)$ of points indexed by the
finite set $A$ in a smooth manifold
$M$ is the subset of $M^A$ given by
\begin{equation}
\oC_A(M) = \{ (x_a)_{a \in A} \in M^A:
	x_a \ne x_b \text{ for }a,b\in A, a \ne b \}
\end{equation}
This is a smooth manifold, but it is not compact if the size of $A$,
$|A|$, is greater than 1.  The
construction is
functorial on the category of smooth manifolds and imbeddings, with,
for $f: N \rightarrow M$ an imbedding and $(x_a)_{a \in A} \in L^A$,
\begin{equation}
\oC_A(f) (x_a)_{a \in A} = (f(x_a))_{a \in A}
\end{equation}
which is in $\oC_A(M)$ since $f$ is an imbedding.
\end{definition}

I sometimes use the notation $\oC_n(M)$ as a shorthand for
$\oC_{\{1,\ldots,n\}}(M)$ when I need to refer to the points
specifically by name.

The compactification of $C_A(M)$ of $\oC_A(M)$ appropriate for this
context was first written down in an
algebro-geometric context\footnote{Using the projective sphere
bundle rather than the sphere bundle so the result is an honest
manifold (or variety).} by Fulton and
Macpherson~\cite{fulton-macpherson:ConfigSpace}, and in the case of
manifolds by
Axelrod and Singer~\cite{Axelrod-Singer:ChernPertII}.

Section~\ref{sec:configdescrip} aims to give an intuitive explanation
and complete description of this wonderful
compactification.  Many details of the construction are deferred
until Section~\ref{sec:configdef}, which is somewhat tangential to the
main thrust of the thesis.  It contains a new definition
of the compactification.
Section~\ref{sec:configvariation}
describes some variations on the compactification necessary in our
context.

\section{A description of the compactification}
\label{sec:configdescrip}

Our basic requirement on the compactification is that the maps
$\phi_{ab}: \oC_A(\RR^3)
\rightarrow S^2$ must extend smoothly to the compactification for all
$a, b\in A$.  Since the compactification should be compact, every path
in $\oC_A(\RR^3)$ must have a limit.\footnote{
Actually, I won't allow paths
that approach spatial infinity for the moment, so ``compactification''
is a slight misnomer.  See Section~\ref{sec:configvariation} for a
description of the true compactification.}
Every such path has a
well-defined limit in $(x_a) \in (\RR^3)^A$.
Group $A$
into classes that map into the same point in $\RR^3$ in the limit.
Then the
direction between elements $a, b$ in different classes is determined
by the limit of $x_a$ and $x_b$; so to guarantee that the
corresponding $\phi_{ab}$ extends continuously, the compactification
will come with a map $C_A(\RR^3) \rightarrow (\RR^3)^A$.

Now onsider a path in $C_A(\RR^3)$ in which the points of just one
subset $B$ of size greater than one approach the same point in
$\RR^3$.  The $\phi_{ab}$ with $a, b \in B$ are completely determined
by the relative configuration of the points in $B$: the positions of
the points in $B$, modulo overall translation and scaling.  This
relative configuration comes with a natural map to the space of points
$(\RR^3)^B$ whose coordinates are not
all identical, modulo translation and scaling; this latter space is
compact (it is the sphere $(((\RR^3)^B/\Delta) \setminus
\{0\})/\RR^+$, where
$\Delta$ is the diagonal in $(\RR^3)^B$), so every path has a limit
point in it.
So to help define the $\phi_{ab}$, those portions of $C_A(\RR^3)$ in which
the points of $B$ approach each other will come with a map to the
space of relative configurations, which I'll call $C_B(\calT\,\RR^3)$.

But the relative configuration might not be enough: the relative
configuration might approach a degenerate relative configuration, in
which some set $C \subset B$ come together.  This corresponds to
all the points in $B$ approaching each other, but the points in $C$
approaching each other faster.  But \emph{all}
the points of $B$ cannot come together in the limit, so we always gain
something: $C$ will be strictly smaller than $B$.  So if the points in
$C$ come together in $C_B(\calT\,\RR^3)$, new coordinates, the
relative configuration of the points in $C$, will help define more of
the $\phi_{ab}$.

In general, the coordinates of $C_A(\RR^3)$ will depend
on the particular point.
The points in some subsets $B_i \subset A$
which approach each other at the top level (only counting maximal $B_i$),
and the directions between points in
different $B_i$ are determined by a macroscopic configuration, a point
in $C_{A/B_i}(\RR^3)$.
($A/B_i$ means $A$ with points in each $B_i$ identified.)  Then, for each
$B_i$, the micrscopic configuration near the $B_i$, the limit point in
$C_{B_i}(\calT\,\RR^3)$ determine some more
directions.  Possibly some $C_{ij} \subset B_i$ still approach each other;
then the submicrscopic configuration at the $C_{ij}$, the limit point
in $C_{C_{ij}}(\calT\,\RR^3)$ determines yet more
directions, and so forth.

The set of subsets $\fS$ that occur in this construction characterize the
coordinates used.  The valid $\fS$ are those
which consist of subsets $S_i$ of $A$ of size $\ge 2$ for which every
distinct pair $S_i, S_j$ satisfy one of
\begin{itemize}
\item $S_i \cap S_j = \emptyset$
\item $S_i \subset S_j$
\item $S_j \subset S_i$
\end{itemize}
These $\fS$ characterize the faces $C_{A,\fS}(\RR^3)$ of $C_A(\RR^3)$,
which is a
\emph{manifold with corners}, a generalization of a manifold with
boundary.  (See below for a precise definition.)  The reader should
verify that the dimension of the face
$C_{A,\fS}(\RR^3)$ is $3|A| - |\fS|$, or codimension $|\fS|$.
($|\fS|$ is the number of distinct subsets of $A$ in $\fS$.)  The
stratum correspond to $\fS$ is, intuitively, the set of limit points
in the open configuration space in which points in each subset $S_i$
approach each other at comparable sppeds.

Although we won't use the compactification for any spaces other than
$\RR^3$ or $S^1$, it's worth pointing out that a similar
compactification can be defined for any $n$-dimensional
manifold $M$.  The
one essential difference is that, without the Euclidean structure of
$\RR^3$, the relative configurations $C_B(\calT M_x)$ are defined only
locally.  Specifically, a coordinate chart $(U, \psi)$ near $x$ yields a
homeomorphism $\psi'$ from $U$ to a neighborhood of $0$ in the tangent
space to $M$ at $x$, $\calT M_x$,
whose Jacobian at $x$ is the identity.  In $\calT M_x$, the limiting
configuration of the points in $B$ is well defined (taking values in
$C_B(\calT M_x) = (((\calT M_x)^B / \Delta) \setminus \{0\})/\RR^+$:
non-constant maps from $B$ 
into the tangent space to $M$ at $x$, modulo overall translation and
scaling).  The limit point in $C_B(\calT M_x)$ of a path in $\oC_A(M)$
in which the points in $B$ all approach the point $x$ does not depend
on the coordinate chart.

\section{Precise definitions}\label{sec:configdef}

It's now time for some proper definitions.  This section is somewhat
off the main track of my thesis, and is somewhat independent of the
other material.  Statements are made in full generality that are not
needed for these knot invariants.  But the results are of some
independent interest.  In particular, I give a new and somewhat simpler
definition of the compactification of configuration space.

As I mentioned, the compactification will be a manifold with corners.

\begin{definition}\label{def:mancorners}
A (smooth) $n$-dimensional \emph{manifold with corners} $M$ is a Hausdorff
topological space
covered by open sets, each homeomorphic to
\begin{equation}
\{x_1,\ldots,x_n) \in \RR^n: x_1 \ge 0, \ldots, x_m \ge 0\} =
	(\bar\RR^+)^m\times\RR^{(n-m)}
\end{equation}
for some integer $0 \le m \le n$.  The transition functions between
two such open sets must extend to a smooth function on an open subset
of $\RR^n$.
\end{definition}

This definition allows us to define smooth functions, smooth
differential forms, the tangent space, etc. The general principle is
that smooth functions (et al) are required to have a smooth extension
to a neighborhood of the coordinate chart.

Any manifold with corners $M$ can be written as the disjoint union of
$m$-dimensional manifolds, for $0 \le m \le n$, consisting of those
points $p\in M$ which have a neighborhood isomorphic to
$(\bar\RR^+)^m\times\RR^{(n-m)}$, with $p$ mapping to the origin.  These
spaces (or the connected components of them) are called the
\emph{strata} of $M$.

The notion of \emph{blowup} will be central in the construction of the
space $C_n(M)$.  Intuitively, the blowup of a manifold  $M$ along a
submanifold $L$, remove $L$ from $M$ and glue in
a replacement that records the directions points approach $L$.  The
replacement is ${\cal S\cal N}L$, meaning the sphere bundle of the
normal bundle of $L$.  The \emph{normal bundle} of $L$ is ${\cal
T}M/{\cal T}L$---that is, tangent vectors of paths approaching $L$, modulo
translation along $L$ itself.  And the \emph{sphere bundle} $\cal S
\cal V$ of a
vector bundle $\cal V$ is, over each point $x$, ${\cal S\cal V}_x = {\cal
V}_x\setminus\{0\} / \RR^+$.  Thus ${\cal S\cal N}L$ removes paths that
sit on $L$ and takes the quotient by the length of the tangent vector in ${\cal
N}L$.

\begin{definition}
The blowup of a smooth manifold with corners $M$ along a closed imbedded
submanifold with corners $L$ is the manifold with boundary $\Bl(M,L)$
that is $M$ with
$L$ replaced by
those points of $\mathcal S(\mathcal N(L))$ that are actually the images of
paths in $M$.  In particular, there is a natural smooth map
$\pi: \Bl(M,L) \rightarrow M$ and a partial inverse
$M \setminus L \rightarrow \Bl(M,L) \setminus \pi^{-1}(L)$.

To construct a manifold diffeomorphic to this blowup explicitly,
pick a Riemannian metric on $M$ and find an $\epsilon > 0$ so that
$B_\epsilon(L) = \{x \in M: d(x, L) < \epsilon\}$ is a tubular
neighborhood of $L$.  Then $M - B_{\epsilon/2}(L)$ is
diffeomorphic to our desired blowup.

If $M$ has a Riemannian structure, then $\mathcal S(\mathcal N(L))$ can be identified
with the unit vectors in $\mathcal TM|_L$ that are perpendicular to $\mathcal TL$.
\end{definition}

Note the similarity with the coordinates on configuration space $C_A(M)$
constructed above.
What I called the ``relative configuration'' of $B$ points at $x \in
M$ is an element on the boundary of the blowup of $M^B$ along the
main diagonal $\Delta = \{(p, \ldots, p) \in M^B, p \in M\}$.
The configuration space $C_A(M)$ should therefore be obtained, in some sense,
by blowing up $M^A$ along each of its diagonals $\Delta_{B}$, $B
\subset A, |B| \ge 2$.  But it's impossible to do all the blowups at
once, and successive blowups are slightly tricky.
If you blow up
along $\Delta_B$ before $\Delta_C$, for $B$ a proper subset of $C$,
then $\Delta_C$
is part of the boundary of the partial blowup.  The blowup along
$\Delta_C$ can be done, but it's easier to just require that
$\Delta_C$ be blown up before $\Delta_B$ is.

With this restriction, when we blow up a diagonal $\Delta_B$ it
will always intersect the interior of the partial blowup, and we can
define successive blowups using the \emph{proper transform}: the
proper transform of $\Delta_B$ will be the closure of $\Delta_B$ minus
its subdiagonals within the partial blowup.  Blowing up
$\Delta_B$ means blowing up the proper transform of $\Delta_B$.

\begin{figure}
\centerline{\input{draws/interline-blowup.tex}}
\caption{The successive blowup along the two intersecting lines in
$\RR^3$ on the left is shown on the right.  Blowups in the other order
would have produced a different result.}
\label{fig:interlineblowup}
\end{figure}
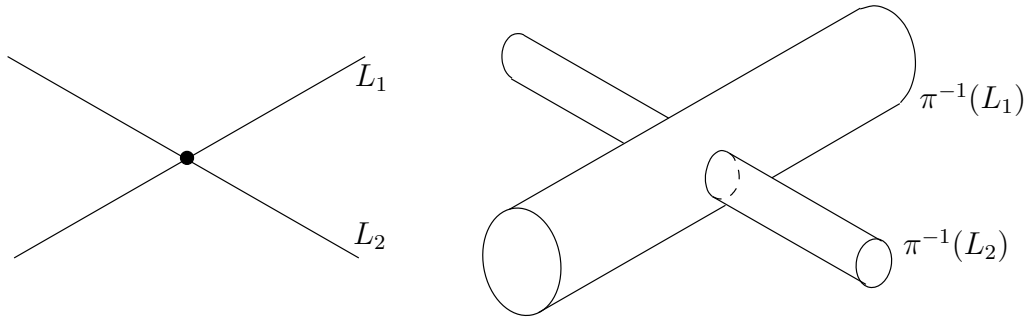

Another worry is that the resulting configuration space might depend
on the order chosen for the blowups, which would mean that it wouldn't
be truly functorial.  For an example of blowups that don't commute,
see Figure~\ref{fig:interlineblowup} for the blowup of two
intersecting lines in $\RR^3$.

Another example that's more relevant
to our situation but harder to visualize is the blowup of three 2-planes
in $\RR^4$ intersecting at 0.  (In fact, this is exactly what the
intersection of the
diagonals $\Delta_{12}, \Delta_{23},$ and $\Delta_{13}$ in $(\RR^2)^3$
looks like, after taking the quotient by overall translation.)
The blowup along the first two
(in either order)
will be isomorphic $\Bl(\RR^2,0) \times \Bl(\RR^2,0)$.  In particular,
the stratum lying over 0 is isomorphic to $S^1 \times S^1$.  But the
third 2-plane will intersect the boundary of the partial blowup only
in the stratum lying over $0$, and in some diagonal on the torus $S^1
\times S^1$.  Blowing up along this plane turns the stratum over 0
into a torus minus a diagonal, which isn't very nice and won't be
symmetric.  

With the restriction that all the $\Delta_C$ for $B \subset C$ and $B
\ne C$ be
blown up before $\Delta_B$, all diagonals become disjoint before you
get around to blowing them up.
Blowups along disjoint submanifolds clearly commute (and can in fact
be done simultaneously), so any permissible order of doing the blowups
would yield the same result.
But this is not true.  For instance, after
blowing up $M^4$ along each of the diagonals $\Delta_S$, $|S| = 3,4$,
the diagonals $\Delta_{12}$ and $\Delta_{34}$ will still
intersect.  But this intersection seems somehow ``nicer'' than the
nasty cases given above.  This can be made precise.

\begin{definition}\label{def:generic}
A collection $N_i$, $1 \le i \le k$ of submanifolds of a manifold
 $M$ intersect \emph{generically} at a point $p\in M^o$ 
if the map
\begin{equation}
\bigoplus_{i=1}^k \mathcal TN_i|_p \rightarrow (\mathcal TM|_p)^k/\Delta
\end{equation}
is surjective, where $\Delta$ is the image of the diagonal imbedding
of $\mathcal TM|_p$ in
$(\mathcal TM|_p)^k$ and the map is the direct sum of the imbeddings of the
$\mathcal TN_i|_p$ in $\mathcal TM|_p$.

(The motivation is that this map is surjective exactly when the
appropriate linear equations can be solved to find a common
intersection point for  small displacements of the $N_i$).

The submanifolds $N_i$ intersect generically if, for every $p \in M$
in at least one of the $N_i$, the $N_i$ such that $p \in N_i$
intersect generically at $p$. 

Equivalently, the manifolds $N_i$ intersect
generically at the point $p\in M$ if and only if $p$ has a neighborhood $U$
such that
\begin{equation}\begin{split} \label{eq:genericprod}
U &\simeq L_0 \times L_1 \times \cdots \times L_k \\
N_i \cap U &= \pi_i^{-1}(p_i),\quad 1\le i\le k
\end{split}\end{equation}
for some  manifolds $L_0, L_i$ and points $p_i \in L_i$,
$0\le i\le k$, so
that $p = p_0 \times p_1 \times\cdots\times p_k$.
($\pi_i$ is the projection of $L \times L_1 \times \cdots\times L_k$
onto the factor $L_i$.)
\end{definition}

\begin{proposition}
If the $N_i$, $1 \le i \le k$, are proper submanifolds of a manifold
with corners $M$ that intersect
generically, then successive blowups along the $N_i$ yield a result
independent
of the order.
\end{proposition}
\begin{proof}
It's enough to proceed locally, in the neighborhood of a point $p \in
M$.  For ease of notation, assume that all the $N_i$ intersect at $p$.
Choose a neighborhood $U$ of $p$ that can be written as in
Definition~\ref{def:generic} as $U \simeq L\times
L_1\times\cdots\times L_k$, with $N_i = \pi_i^{-1}(p_i)$.  Then the
blowup of $U$ along $N_1$  is $L \times \Bl(L_1,p_1)\times L_2
\times\cdots\times L_k$, the blowup of this space along $L_2$ is
$L\times\Bl(L_1,p_1)\times\Bl(L_2,p_2)\times\cdots\times L_k$,
and so forth; after all the blowups, the result is
$L\times\Bl(L_1,p_1)\times\cdots\times \Bl(L_k,p_k)$,
independently of the order they were preformed.

(The condition that the $N_i$ be proper submanifolds is only required
because we can't blow up a manifold $M$ along $M$ itself.)
\end{proof}

So the notion of generic intersection seems to correspond to intuition
about ``niceness'' of the intersection.  We're ready to define the
space $C_A(M)$, in stages from $M^A$.  At each stage, we've blown up
some collection $\Delta_{B_i}$ of diagonals of $M^A$.  The next stage
is the blow up of any diagonal $\Delta_C$ of $M^A$ for which all
proper supersets of $C$ are contained in the $B_i$.  To prove the
result doesn't depend on the choice of which $\Delta_C$ to blow up
at each stage, it suffices to check that
adjacent steps commute.

\begin{proposition}
If the diagonals corresponding to all proper supersets of $C$ and $D$
have been blown up at some stage in the construction of $C_A(M)$, but
not $\Delta_C$ and $\Delta_D$, then $\Delta_C$ and $\Delta_D$
intersect generically.
\end{proposition}
\begin{proof}
If $C \cap D \ne \emptyset$, then their intersection at this stage
will be empty.  In this case, $\Delta_C \cap \Delta_D = \Delta_{C \cup
D}$, which, by hypothesis, we've already blown up; and in the
coordinates over $\Delta_{C \cup D}$, there exist points with either
all the points in $C$ identical, or all the points in $D$ identical,
but not both.

On the other hand, if $C \cap D = \emptyset$, $\Delta_C$ and
$\Delta_D$ will intersect, but the intersection will be generic.

at a given stratum, want to reduce to each subst. individually.
Right: enough to show that each tangent factor is a product as we want.
\end{proof}

For any embedding $f:
N \rightarrow M$, the diagonals $\Delta_B$ of $N^A$ are the
intersections of the diagonals $\Delta_B$ of $M^A$ with $N^A$, with a
generic intersection.  Since the sequence of blowups is the same in
the two cases, the blowup will be functorial.

\section{Further variations}\label{sec:configvariation}

The configuration spaces of interest to us differ from the
spaces $C_A(\RR^3)$ defined above in two
ways.  First, as previously mentioned, $C_A(\RR^3)$ is not actually compact
as $\RR^3$ is not compact.
Instead of $C_A(\RR^3)$, the proper compactification is
$C_{A \cup\{\infty\}}(S^3)$ for some (fixed) point $\infty \in S^3$.
$\RR^3 \simeq S^3 \setminus \infty$, so $\oC_A(\RR^3) \simeq \oC_{A
\cup \{\infty\}}(S^3)$, and $C_{A \cup\{\infty\}}(S^3)$ is a proper
compactification of $\oC_A(\RR^3)$.  The maps $\phi_{ab}$ extend
continuously to this compactification as well; the direction a point
$a \in A$
approaches $\infty$ determines the direction from $a$ to points that do not
approach $\infty$, and the relative configuration at $\infty$
determines directions between points that approach $\infty$.

Since points approaching $\infty$ turn out to be mostly irrelevant for
the purposes of this thesis, in the future I'll write $C_A(\RR^3)$ to
mean $C_{A\cup\{\infty\}}(S^3)$,
leaving the strata at $\infty$ implicit.

The other difference is the presence of the knot.  As explained
above, this is a functorial construction, so a knot $K: S^1
\rightarrow \RR^3$ yields an imbedding $C_A(K)$ (of
manifolds with corners) of
$C_A(S^1)$ in $C_A(\RR^3)$, which defines
$\phi_{ab}$ for
points on the knot and use these to construct certain integrals.  But
as explained in Section~\ref{sec:introduction}, for any but the
simplest integrals this is not
enough to construct an invariant.  Most invariants involve integration
of certain points over $\RR^3$ and other points over $S^1$.  To do
this properly, we need the notion of a \emph{restricted configuration
space} $\oC_{A;B}(\RR^3; K)$ of maps from $A$ to distinct points in
$\RR^3$, with the
points in $B$ landing in the submanifold $K(S^1)$.  (Note that points
in $A \setminus B$ range over the complement of the $x_b$, $b \in B$,
not the complement of the knot $K(S^1)$.)
Its closure in
$C_A(\RR^3)$, $C_{A;B}(\RR^3;K)$, is the desired compact space.

Since the knot $K$ is one-dimensional, the spaces $C_{A;B}(\RR^3;K)$ as
I've defined them are not connected for $|B| > 2$.  The (cyclic) order the
points in $B$ occur around $K$ is significant.  In this notation, $B$
will therefore be taken to include a cyclic order of its points.

It's important to note that $C_{A;B}(\RR^3;K)$ is smoothly fibered over the
space of knots $\knots$.  $\oC_{A;B}(\RR^3;K)$ is a subset of
$(\RR^3)^{A\setminus B} \times (S^1)^B$ defined by equations like $x_a
\ne x_b$ and $x_a \ne K(s_b)$ and so is fibered over $\knots$.  Since
$C_{A;B}(\RR^3;K)$ is just the closure of $\oC_{A;B}(\RR^3;K)$ in
$C_A(\RR^3)$, it too is fibered over $\knots$.  The notation $C_{A;B}$
without qualifications will mean this bundle over $\knots$.

Finally, since the codimension 1 strata of $C_{A;B}$ become
important in Stokes' theorem, a brief description is
appropriate.  As usual, the strata are parametrized by subsets $D
\subset A$.  If $D \subset A \setminus B$ consists of only free
vertices, then the stratum can be written $C_{(A;B),D} = C_{(A/D);B}
\times C_D(\calT\,\RR^3)$: that is, this stratum of configuration
space is the product of the macroscopic configuration with the
microscopic configuration.

On the other hand, if $B \cap D \ne \emptyset$ (so $D$
contains some points on the knot, necessarily a contiguous subset in
the cyclic order on $B$),
$C_{(A;B),D}(\RR^3;K)$ is fibered over $C_{(A/D);(B/D)}(\RR^3; K)$, with
fiber $C_{D;(B\cap D)}(\calT\,\RR^3; \dot K(x_D))$.
This is again the macroscopic configuration and the microscopic
configuration, but now the microscopic configuration is relative
the line $\dot K(x_D)$, the derivative of $K$ at the common limit of $D$.
This latter
space is the intersection of the closure of $\oC_{D;(B\cap
D)}(\RR^3;K)$ with the stratum $C_{D,D}(\RR^3) \simeq
C_D(\calT\,\RR^3)$ near the point $x_D$ 
(the common limit point of the points in $D$).
Specifically, this is the space of non-constant maps $D \rightarrow \RR^3$
with points in $D \cap B$ in the subspace spanned by $\dot
K(x_D)$, modulo translations in the direction $\dot K(x_D)$ (so that
points in $D \cap B$ remain in the right subspace) and overall
scaling.

%% file: draws/interline-blowup.tex
\begin{picture}(0,0)%
\epsfig{file=draws/interline-blowup.pstex}%
\end{picture}%
\setlength{\unitlength}{3947sp}%
\begingroup\makeatletter\ifx\SetFigFont\undefined
\def\x#1#2#3#4#5#6#7\relax{\def\x{#1#2#3#4#5#6}}%
\expandafter\x\fmtname xxxxxx\relax \def\y{splain}%
\ifx\x\y   
\gdef\SetFigFont#1#2#3{%
  \ifnum #1<17\tiny\else \ifnum #1<20\small\else
  \ifnum #1<24\normalsize\else \ifnum #1<29\large\else
  \ifnum #1<34\Large\else \ifnum #1<41\LARGE\else
     \huge\fi\fi\fi\fi\fi\fi
  \csname #3\endcsname}%
\else
\gdef\SetFigFont#1#2#3{\begingroup
  \count@#1\relax \ifnum 25<\count@\count@25\fi
  \def\x{\endgroup\@setsize\SetFigFont{#2pt}}%
  \expandafter\x
    \csname \romannumeral\the\count@ pt\expandafter\endcsname
    \csname @\romannumeral\the\count@ pt\endcsname
  \csname #3\endcsname}%
\fi
\fi\endgroup
\begin{picture}(5749,1966)(5509,-4596)
\put(7700,-4124){\makebox(0,0)[lb]{\smash{\SetFigFont{12}{14.4}{rm}$L_2$}}}
\put(11258,-3271){\makebox(0,0)[lb]{\smash{\SetFigFont{12}{14.4}{rm}$\pi^{-1}(L_1)$}}}
\put(11153,-4187){\makebox(0,0)[lb]{\smash{\SetFigFont{12}{14.4}{rm}$\pi^{-1}(L_2)$}}}
\put(7710,-3134){\makebox(0,0)[lb]{\smash{\SetFigFont{12}{14.4}{rm}$L_1$}}}
\end{picture}

%% file: tying.tex
\chapter{Tying it Together}
\label{sec:tying}

In this chapter I bring together a number of loose ends, completing
the proof that the invariants I'll construct will be Vassiliev
(Section~\ref{sec:integralvas}), explaining why the variation of an
integral is the integral restricted to the boundary and how the same
cancellations can prove that a number of different integrals are
invariants and give the same value (Section~\ref{sec:pushforward}),
making some
notes on the form of the integrals resulting
(Section~\ref{sec:integrals}), 
completing the proof that the contributions from principal faces
cancel for an integral $I(w)$ corresponding to a weight diagram $w$
(Section~\ref{sec:cancel}), and introducing \emph{tinkertoy diagrams}
(Section~\ref{sec:tinkertoy}).

\section{Integrals are Vassiliev}
\label{sec:integralvas}

We can define the configuration space $C_{A;B}(\RR^3;K)$ of points on
an immersed curve $K$ with $k$ double points: it is the closure in
$C_A(\RR^3)$ of $\oC_{A;B}(\RR^3;K)$, the space of $B$ distinct points on
the immersed curve and $A\setminus B$ points
in space. This is a perfectly fine compact
manifold with corners, so configuration space integrals make
sense for curves with double points as well.

Given any configuration space integral on a curve $K$ with $k$ double
points, suppose we approach a curve $K^x$ with $k+1$ double points by
varying $K$.
As the strands near $s_1, s_2 \in S^1$ approach each other (with the
tangents at $s_1$ and $s_2$, $\dot
K(s_1)$ and $\dot K(s_2)$, remaining fixed),
the submanifold $C_{A;B}(\RR^3;K)$ of $C_A(\RR^3)$ approaches a
limiting position.  The limiting position is not a smooth submanifold
of $C_A(\RR^3)$.
But there is an open dense subset of the limiting position that is
smooth, and so the limit of
the integral will be equal to the integral over this subset of the
limiting position.

There are a number of connected components to this
smooth portion of the limiting position, indexed by the set $D \subset
A$ of points that get very close to the crossing, together with the
subsets $D_1 \subset D \cap B$ and $D_2 \subset D \cap B$ that
approach the crossing along the strands near $s_1$ and $s_2$
respectively.  (Of necessity, $D_1$ and $D_2$ will each be connected
subsets of $B$ and $D \cap B = D_1 \amalg D_2$.)  This portion of the
limiting position then lies in the stratum $C_{A,D}(\RR^3)$ of
$C_A(\RR^3)$ in
which the points $D$ approach each other (at comparable speeds) and is
the product of
the macroscopic configuration space $C_{A/D; B/D}(\RR^3; K^x)$ with
a microscopic configuration space, whose details don't matter much,
but is
the image of
[configurations of points $D_1$ on a line in the direction $\dot
K(s_1)$, points $D_2$ on a line in the direction $\dot K(s_2)$,
and points $D \setminus (D_1 \cup D_2)$ in $\RR^3$] in
$C_D(\calT\,\RR^3)$.

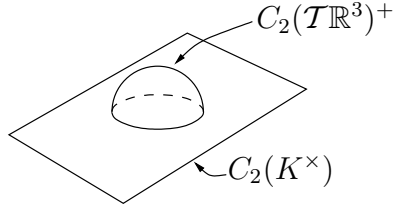
\begin{figure}
\centerline{\input{draws/config-crossing.tex}}
\caption{Local situation of two points near the limit of a crossing,
drawn as a subset of 
$C_2(\RR^3)$ modulo overall translation.  The flat section represents
points macroscopically separated and the hemisphere in the center is the
portion of $C_2(\calT\,\RR^3)$ that is covered by the limit.}
\label{fig:configcrossing}
\end{figure}

 See Figure~\ref{fig:configcrossing} for a picture
of the situation in the most relevant case, when $D$ has two points,
one on each strand.

So the limiting integral can be broken into components with different
subsets $D$ very near the crossing.  But the integral over a component
with 0 or 1 points near the crossing will clearly not depend on the
direction of the crossing.  Hence the $k+1$st derivative of a knot
invariant with $2k$ crossings 
will be 0.  This completes the proof of
Proposition~\ref{prop:integralvas}.

\section{On the Pushforward}
\label{sec:pushforward}
Our integrals are of the general form
\begin{equation}
\int_{C_{A;B}(\RR^3;K)} \omega
\end{equation}
where $\omega$ is an $n$-form on $C_{A;B}$,
with $n = 3|A| - 2|B| = \dim C_{A;B}(\RR^3;K)$.  (Recall that
$C_{A;B}$ is the fiber bundle over the space $\knots$ of knots whose
fiber over $K$ is $C_{A;B}(\RR^3;K)$.  For us, $\omega$ will be the
pullback of a top-dimensional form on some manifold $M$, and so will
be closed.)
This type of expression is called \emph{integrating over the fiber} or
the \emph{pushforward}.

\begin{definition}
Given any $n+k$-form $\omega$ 
 on a fiber bundle $\pi: B
\rightarrow X$ whose fibers are oriented $n$-dimensional manifolds
with corners, the \emph{pushforward} of $\omega$ along $\pi$, denote
$\pi_*\omega$, is
the $k$-form on $X$ whose
value on a $k$-chain $c$ is
\begin{equation}
\int_c \pi_*\omega = \int_{\pi^{-1}(c)} \omega
\end{equation}
\end{definition}

When $B$ is a smooth manifold (without boundary) the pushforward of a closed
form will be closed.  But
when $B$ is a manifold with corners, this is not true; instead, we
have

\begin{proposition}\label{prop:dpushforward}
For bundle $\pi: B\rightarrow X$ whose fibers are smooth oriented
manifolds with
corners with boundary $(\partial \pi): \partial B \rightarrow X$
(with the orientation induced from $B$),
\begin{equation}
d\pi_*\omega = \pi_*d\omega - (\partial \pi)_*\omega
\end{equation}
\end{proposition}
\begin{proof}
For any $k+1$-chain $c$ in $X$, we have
\begin{equation}
\begin{split}
\int_c \pi_*d\omega &= \int_{\pi^{-1}(c)} d\omega\\
	&= \int_{\partial(\pi^{-1}(c))} \omega\\
	&= (\int_{\pi^{-1}(\partial c)} +
		\int_{(\partial\pi)^{-1}(c)}) \omega\\
	&= d\pi_*\omega + (\partial\pi)_*\omega
\end{split}
\end{equation}
using Stokes' theorem and $\partial(\pi^{-1}(c)) = \pi^{-1}(\partial
c) \amalg (\partial\pi)^{-1}(c)$.
\end{proof}

This proposition justifies my claim that the variation of an
integral is (up to sign) the integral restricted to the boundary:
since the forms $\omega$ 
we integrate are closed\footnote{
Even though the dimension of the fiber is equal to the degree of
$\omega$, this is not a trivial statement, since it requires $\omega$
to be closed as a form on the bundle.} (being the pullback of a closed
form), $\pi_*d\omega$ will be 0, and $d\pi_*\omega$, the variation of
the integral of $\omega$, will be $-(\partial\pi)_*\omega$.

This proposition also tells us what happens as we vary the $n$-form
$\omega$.  If we pick a form $\omega'$ in the same cohomology class as
$\omega$, then $\omega' - \omega = d\alpha$ for some $n-1$-form
$\alpha$.  Applying the above proposition, $\pi_*\omega' - \pi_*\omega
= \pi_*d\alpha = d\pi_*\alpha + (\partial\pi)_*\alpha$; but
$\pi_*\alpha$ is a $-1$-form on $X$ and hence 0.  Thus the difference
between the integrals of $\omega$ and $\omega'$ will be
$(\partial\pi)_*\alpha$ for some form $\alpha$.

Therefore an integral will both
\begin{itemize}
\item be an invariant under deformation of
the knot and 
\item be independent of the form $\omega$ within a cohomology
class
\end{itemize}
if the boundary of the bundle $B$ is 0 in some sense.  More
precisely, since the forms we use are all $\phi^*\psi$ for some
$n$-form $\psi$
on a manifold $M$ (which will be a product of 2-spheres)
the image of the boundary of $B$ in this product should be 0 (as an
$n-1$-chain).  But
since there is not yet any sensible homology theory with general
manifolds with corners, I'll state this in a different form.

Precisely, the arguments I'll use to show the vanishing or
cancellation of the pushforward of $\phi^*\omega$ along a given
non-anomalous face
$\partial B$ will
either be
\begin{itemize}
\item degeneration: a proof that the image of each fiber of $\partial
B$ in $M$ lies 
in an $n-2$ dimensional submanifold; or
\item cancellation: a map $\Phi$ from $\partial B$ to some other
boundary component $\partial
B'$ that is orientation reversing but commutes with the maps to $M$
(for example, or else $\Phi$ might be orientation preserving but the
two faces
appear with an opposite numerical factor, or 3 (or more) faces might
cancel).
\end{itemize}
In either case, it follows that both desiderata above hold.  If the
face degenerates, then $(\partial\pi)_*\phi^*\psi$ evaluated on any
1-chain is zero, as the image of the fiber will be a
codimension 1 submanifold; and $(\partial\pi)_*\phi^*\alpha$, for an $n-1$
from $\alpha$ on $M$, will be zero as a function on $\knots$ since the image
of the fiber is of dimension $n-2$.  The arguments are even easier in
the case of cancellation.  

Separate arguments are needed for some of the anomalous faces.  These
are discussed in
Chapter~\ref{sec:anomaly}.

\section{The Form of the Integrals}
\label{sec:integrals}

As mentioned in Chapter~\ref{sec:introduction}, the configuration
space integrals are the integrals of the pullback of the canonical form on
$(S^2)^e$ via the map $\phi(\Gamma): C_{V;c} \rightarrow (S^2)^e$
defined by some oriented trivalent graph $\Gamma$.  Recall that the
canonical form on $S^2$ is
\begin{equation}
\begin{align}
\omega 
       &= \frac{\epsilon_{\mu\nu\sigma}}{8\pi}\,
		\frac{x^\mu\,dx^\nu\,dx^\sigma}{|x|^3}
\end{align}
\end{equation}
Pulling this back to $C_{A;B}$ via the map $\phi_{ab}: C_{A;B}
\rightarrow S^2$, the direction map between $x_a$ and $x_b$, yields
\begin{equation}\label{eq:theta_ab}
\theta_{ab} = \phi_{ab}^* \theta =
	\frac{\epsilon_{\mu\nu\sigma}}{4\pi}\,\frac{x_b-x_a}{|x_b-x_a|^3}
		\,(\frac{1}{2}dx_a^\nu\,dx_a^\sigma
		   -dx_a^\nu\,dx_b^\sigma
		   +\frac{1}{2}dx_b^\nu\,dx_b^\sigma)
\end{equation}

\begin{figure}
\centerline{\input{draws/bad-graph.tex}}
\caption{A bad graph.  It is not trivalent, yet yields a
non-degenerate map from configuration space to the product of
spheres.}
\label{fig:badgraph}
\end{figure}
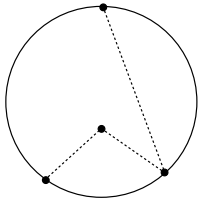

All the invariants considered in this thesis have trivalent vertices,
both internal and on the knot (counting the knot as 2 edges incident
to a vertex
on the knot).  A simple dimension check verifies that these yield
functions on knots. (The dimension of the configuration space is
(\# vertices on knot) + 3(\# vertices off knot), while the dimension of
the product of spheres is 2(\# propogators); since each
propogator has two
endpoints while each vertex is incident to 1 or 3 propogators, these match.)
There are other graphs that yield non-zero functions on the space of
knots.
See figure~\ref{fig:badgraph} for
an example.  These graphs will not be considered in this thesis.
But since every configuration space integral is a
Vassiliev invariant (Proposition~\ref{prop:integralvas}) and I give an
explicit construction of an integral for every Vassiliev invariant, no
new invariants can be constructed from them.

For a trivalent diagram and the canonical form on the product of
$S^2$'s, the integral can be written in another way.  The product of the
various $\theta_{ab}$'s as in~\ref{eq:theta_ab} is a sum of terms with
0, 1, or 2 $dx_a$'s
and a complementary number of $dx_b$'s.  But note that, for the
integral of this form over $C_{A;B}(\RR^3; K)$ not to vanish, exactly
one $dx_a$ must appear for each $a \in B$
on the knot and exactly three $dx_a$'s must appear for each $a \in A
\setminus B$ in $\RR^3$.  A little combinatorial juggling
shows that only
the terms of the various $\theta_{ab}$ with exactly one $dx_a$ and one $dx_b$
give a non-zero contribution to the integral.  Note that this is only
true for the integral over the fiber. In particular, the term you
integrate here is not a closed form on the bundle, so computing the
variation with this expression involves more than just restriction to
the boundary.

For example, one of
the integrals in the degree 2 Vassiliev invariant is (from the
introduction) $\int_{\oC_{4;3}} \theta_{14}\theta_{24}\theta_{34}$.
By the above remarks, this integral can be rewritten
to look like
\begin{multline}
\label{eq:canonintegral}
\int ds_1\,\frac{dx_1^{\mu_1}}{ds_1}\,
\int ds_2\,\frac{dx_2^{\mu_1}}{ds_2}\,
\int ds_3\,\frac{dx_3^{\mu_1}}{ds_3}\,
\int d^3 x_4\,\epsilon^{\nu_1\nu_2\nu_3}\\
\Delta_{\mu_1\nu_1}(x_4-x_1)
\Delta_{\mu_2\nu_2}(x_4-x_2)
\Delta_{\mu_3\nu_3}(x_4-x_3)
\end{multline}
where $\Delta_{\mu\nu}(x)$ is the propogator written as a tensor:
\begin{equation}
\Delta_{\mu\nu}(x) =
\frac{\epsilon_{\mu\nu\sigma}}{4\pi} \frac{x^\sigma}{|x|^3}
\end{equation}

But note what's required for this type of expression.  
In general, assign a tensor index to each propogator leaving each
vertex and
write down $\int ds_a\,\frac{dx_a^\mu}{ds_a}$ for each
point $a$ on the knot, $\int d^3 x_b$ for each point $b$ in space,
$\Delta_{\mu\nu}(x_b - x_a)$ for each propogator, and a term
$\epsilon^{\mu\nu\sigma}$ for each internal vertex, with the
appropriate choices of tensor index.
As each integral is paired with a corresponding differential form, the order
of the integrals (which also determines the orientation on the space)
is irrelevant.  The only sign choice in an expression like this is the
internal vertex $\epsilon^{\mu\nu\sigma}$ tensors. But this is exactly
an ordering on the edges meeting each internal vertex, which is
\emph{exactly} what's needed to define
weight diagrams.  This gives an algorithm for relating the 

For a more complete discussion of these signs, see
Appendix~\ref{sec:graph}.

\section{Principal Faces Cancel}
\label{sec:cancel}

There are several fixes that need to be applied to the naive
definition of $\phi(\Gamma)$ in the introduction in order to get a
meaningful formula.

To apply the results of Section~\ref{sec:pushforward} on methods of
proving invariance and non-dependence on the form $\omega$ on $M$,
you need maps
from the bundle $C_{A;B}$ to a target manifold that is the same for
all maps.  But the number of propogators in the graphs associated to a
weight diagram varies.  (The number of propogators is the degree plus
the number of internal vertices.)  To fix this, $\phi(\Gamma)$ will be
modified slightly: it will be a map $C_{A;B} \times (S^2)^x
\rightarrow (S^2)^{a+x}$, the product of the naive $\phi(\Gamma)$ with
the identity map on enough spheres so that all these maps have the
same target space.

Equation~\ref{eq:mainthm}, the definition of the integrals considered,
was not actually well defined, since in the definition of the integral
associated to a graph in the Introduction I needed an ordering of the
vertices (to define the orientation on configuration space) and an
ordering of the edges and a direction on each edge (to define the map
to $(S^2)^e$).

The ordering on the edges is actually irrelevant for
the canonical integrals, as the product of canonical forms on each
$S^2$ is symmetric under permutation of the edges.  But for 
general forms that are not products of the same form on each $S^2$
factor, this won't be true. As explained in the next
section, we'll need these more general forms.  To solve this,
$I(\Gamma)$, the integral associated to the graph $\Gamma$,
will include averaging over all orderings of the edges.
The canonical form is also antisymmetric under the antipodal map $N$
on $S^2$; this will be needed below.
Since not all forms on $(S^2)^{e+x}$ will have this property, $I(\Gamma)$
will also be antisymmetrized under the antipodal map on each edge.

(Alternatively, you could view this as a restriction to forms on
$(S^2)^{e+x}$ that are symmetric under permutations of the factors and
antisymmetric under the antipodal map on each sphere.)

To define $I(\Gamma)$ now, you need an ordering on the vertices and a
choice of direction on each edge.  A transposition of two vertices and
changing the direction on an edge both change the sign of the
resulting integral. This defines a notion of orientation on the graph.
But in equation~\ref{eq:mainthm}, the orientation on a graph used was
an orientation on each trivalent edge.  For trivalent diagrams, these
two notions of orientation turn out to be equivalent.  See
Appendix~\ref{sec:graph}, or alternatively take the canonical integral
and rewrite it as described in the last section to get it into a form
as in equation~\ref{eq:canonintegral}, which only depends on an
orientation at each vertex.

Although $I(\Gamma)$ as it stands now is an invariant in the desired
sense, the automorphism factors $|\mathrm{Aut}(\Gamma)|$ in
equation~\ref{eq:mainthm} still
appear a bit mysterious.  To make them transparent, note that the
number of ways of labelling the vertices of an $n$-vertex graph
$\Gamma$ with $\{1, \ldots, n\}$
 is $\frac{n!}{|\mathrm{Aut}(\Gamma)|}$.  So instead of the
equation~\ref{eq:mainthm} as it stands (with its annoying reference to
``isomorphism classes of graphs''), we can equivalently write
\begin{equation}
I(w) = \frac{1}{n!}\sum_{\substack{\mathrm{labelled}\\
	\mathrm{graphs\ }\Gamma}}
{w(\Gamma) I(\Gamma)}
+ \mathrm{(anomaly\ correction\ term)}
\end{equation}

With all these modification, $I(\Gamma)$ is now properly defined.  The
proof that the principal faces cancel now falls into place easily.

First,
the principal face corresponding to the collision of two points not
connected by a propogator or a section of the knot is degenerate: the
stratum of configuration space (for the collision of $a$ and $b$) is
isomorphic to $C_{A/\{a,b\};B/\{a,b\}}(\RR^3, K)
\times C_{\{a,b\}}(\calT\,\RR^3)$, and the map $\phi(\Gamma)$ on this
configuration space is independent of the second component.  Therefore
the image of this stratum must have codimension at least 3 in
$(S^2)^{e+x}$.

The remaining cases will correspond to either the $STU$ or $IHX$
relation.  When two points connected by a propogator collapse, the
stratum of configuration space is again isomorphic to
$C_{A/\{a,b\};B/\{a,b\}}(\RR^3, K)
\times C_{\{a,b\}}(\calT\,\RR^3)$; the map $\phi(\Gamma)$ will be a
product of $\phi(\Gamma_e)$ of the map (corresponding to the labelled
graph $\Gamma_e$ with the edge connecting $a$ and $b$ contracted) with
the identity map on the second component, which is isomorphic to
$S^2$.  When two points connected by a segment of the knot collapse,
the stratum is $C_{A/\{a,b\};B/\{a,b\}}(\RR^3, K) \times
C_{\{a,b\};\{a,b\}}(\calT\,\RR^3; \dot K(s_{\{a,b\}}))$.  This latter
space is just a point and so can be omitted.
$\phi(\Gamma)$ on this stratum is in this case just
$\phi(\Gamma_e)$.  

So in any case, $\phi(\Gamma)$ restricted to this
stratum is $\phi(\Gamma_e)$ cross the identity on some number of
spheres.  Each \emph{labelled} graph $\Gamma_e$ that arises in this
way comes from exactly six
labelled trivalent graphs $\Gamma$ by contraction.  (The vertex
resulting from the contraction will be labelled with the set of the
two contracted vertices. There are 6 ways of partitioning the 4
edges incident to this vertex into 2 (labelled) partitions of two
edges each.)
Depending on whether the contracted vertex is internal or on the knot,
these 6 diagrams correspond to the 3 diagrams in the $IHX$ or $STU$
relation in a natural way.  A simple sign check then verifies that
these diagrams cancel each other out, using the corresponding equation
for the weight diagram.

\section{Tinkertoy Diagrams}
\label{sec:tinkertoy}

In this section I discuss briefly a combinatorial formula for these
configuration space integrals.  This will be expanded in a future
paper.  Right now the main consequence is

\begin{theorem}
For a weight diagram with rational coefficients, the corresponding
configuration space integral invariant is rational.
\end{theorem}

(In fact, it's easy to bound the demoninator by counting everything we
averaged over.)

As I mentioned in Section~\ref{sec:pushforward}, my proof of
invariance of the configuration space integrals proves, at the same
time, that any other form on
$(S^2)^{e+x}$ yields the same invariant.
A
particular choice of the form on $(S^2)^{e+x}$ will yield a counting
formula for this invariant by the following general remarks,
which apply to any invariant for which ``boundaries cancel'' as in
Section~\ref{sec:pushforward}.

Suppose $p\in M$ is a regular value of each
of the maps $f_i:B_i\rightarrow M$, so that each $f_i$ is a
diffeomorphism in a neighborhood of each inverse image of $p$.  (We can
always find such a $p$ by Sard's theorem.)
Then $p$ will
have only a finite number of inverse images (since each of the
$B_i$ is compact) and $f_i$ will
necessarily be non-degenerate (and hence a diffeomorphism) in a
neighborhood $U_{ij}$ of each of the
inverse images $x_{ij}$ of $p$ in $B_i$.  Pick an $n$-form $\omega$ so
that $\int_M \omega = 1$
and $\supp \omega \subset \bigcap_i f_i(U_i)$.  ($\supp \omega$, the
support of $\omega$, is the closure of the set of points in $M$ where
$\omega$ is non-zero.)  Then $\int f_i^* \omega$ is $\sgn f_{i*}(x_{ij})$ near
$x_{ij}$ of $p$, and $\sum_i \int_{B_i} f_i^*\omega$ is
the sum of $\sgn f_{i*}(x_{ij})$ over all inverse images of $p$.
($\sgn f_{i*}(x_{ij})$ is $\pm 1$ depending
on whether $f_i$ is orientation
preserving or reversing near that inverse image.)

Thus
counting the number of inverse
images of a generic point $p\in (S^2)^{e+x}$ (with appropriate signs) yields
a purely
combinatorial formula for each invariant.  These inverse images can be
interpreted as \emph{tinkertoy diagrams}: some nodes, either anchored
to the knot or free in space, and rods between certain pairs of nodes.
The rods are constrained to lie in certain directions, each corresponding
to the projection of $p$ onto the $S^2$ factor for that
edge.  Because of the symmetrization on the $S^2$ factors described in
the previous section, only the \emph{set} of allowable directions
matters; each rod must lie in an allowable direction (or its
negative), with no
direction used more than once.

The sign associated with each tinkertoy diagram can be computed using
the pullback of the
canonical form as in equation~\ref{eq:canonintegral}.  By definition,
the form $ds_1\,ds_2\,ds_3\,d^3x_4$ (for
equation~\ref{eq:canonintegral}, or in general the product of the
$d^3x_a$ and $ds_b$ for $a \in A \setminus B$ and $b \in B$ in the
proper order) is positive (has postive value when integrated over a
small region in $C_{A;B}$), so the total form will be positive exactly
when the numerical factor is positive.  This numerical factor can be
interpreted in
various ways, using the interpretation of
$\epsilon_{\mu\nu\sigma}x^\nu y^\sigma$ as the cross product of the
vectors $x$ and $y$.

For the degree 2 Vassiliev invariant and an almost planar knot (as in
equation~\ref{eq:v2}), this can be simplified somewhat.  The
two types of tinkertoy structures in this case are
\begin{itemize}
\item have two rods
with both endpoints anchored to the knot (and a constraint on the
cyclic order the endpoints occur around the knot) or
\item have a tripod of rods anchored at one end and joined at an internal node.
\end{itemize}
For an almost planar knot, rods with both endpoints anchored to the
knot can only occur near a crossing, while tripods occur exactly when
there is a triangle of points on the knot similar to a given triangle.
So this invariant can be reduced to a sum over certain pairs of
crossings, certain crossings (depending on the directions of the
strands at the crossing), and certain triangles of points on the knot.

By choosing the 3 directions all near a plane, two things happen: the
contributions from single crossings drop out and the triangle degenerates
into counting critical points and points in the same horizontal slice
through the knot with a given ratio.

Dror Bar-Natan has suggested that this
approach can be used to construct a quasi-triangular quasi-Hopf
algebra \textit{\'a la} Drinfel'd~\cite{drin:quasitri,drin:quasihopf}
and in particular to show that these integrals are the same as
those arising from the Knizhnik-Zamoldchikov
connection~\cite[Section~4]{bn:vassiliev}.  This program will be
carried out in a future paper.

%% file: draws/config-crossing.tex
\begin{picture}(0,0)%
\epsfig{file=draws/config-crossing.pstex}%
\end{picture}%
\setlength{\unitlength}{3947sp}%
\begingroup\makeatletter\ifx\SetFigFont\undefined
\def\x#1#2#3#4#5#6#7\relax{\def\x{#1#2#3#4#5#6}}%
\expandafter\x\fmtname xxxxxx\relax \def\y{splain}%
\ifx\x\y   
\gdef\SetFigFont#1#2#3{%
  \ifnum #1<17\tiny\else \ifnum #1<20\small\else
  \ifnum #1<24\normalsize\else \ifnum #1<29\large\else
  \ifnum #1<34\Large\else \ifnum #1<41\LARGE\else
     \huge\fi\fi\fi\fi\fi\fi
  \csname #3\endcsname}%
\else
\gdef\SetFigFont#1#2#3{\begingroup
  \count@#1\relax \ifnum 25<\count@\count@25\fi
  \def\x{\endgroup\@setsize\SetFigFont{#2pt}}%
  \expandafter\x
    \csname \romannumeral\the\count@ pt\expandafter\endcsname
    \csname @\romannumeral\the\count@ pt\endcsname
  \csname #3\endcsname}%
\fi
\fi\endgroup
\begin{picture}(1892,1305)(9722,-2983)
\put(11113,-2806){\makebox(0,0)[lb]{\smash{\SetFigFont{12}{14.4}{rm}$C_2(K^\times)$}}}
\put(11292,-1858){\makebox(0,0)[lb]{\smash{\SetFigFont{12}{14.4}{rm}$C_2(\calT\RR^3)^+$}}}
\end{picture}

%% file: draws/bad-graph.tex
\begin{picture}(0,0)%
\epsfig{file=draws/bad-graph.pstex}%
\end{picture}%
\setlength{\unitlength}{3947sp}%
\begingroup\makeatletter\ifx\SetFigFont\undefined
\def\x#1#2#3#4#5#6#7\relax{\def\x{#1#2#3#4#5#6}}%
\expandafter\x\fmtname xxxxxx\relax \def\y{splain}%
\ifx\x\y   
\gdef\SetFigFont#1#2#3{%
  \ifnum #1<17\tiny\else \ifnum #1<20\small\else
  \ifnum #1<24\normalsize\else \ifnum #1<29\large\else
  \ifnum #1<34\Large\else \ifnum #1<41\LARGE\else
     \huge\fi\fi\fi\fi\fi\fi
  \csname #3\endcsname}%
\else
\gdef\SetFigFont#1#2#3{\begingroup
  \count@#1\relax \ifnum 25<\count@\count@25\fi
  \def\x{\endgroup\@setsize\SetFigFont{#2pt}}%
  \expandafter\x
    \csname \romannumeral\the\count@ pt\expandafter\endcsname
    \csname @\romannumeral\the\count@ pt\endcsname
  \csname #3\endcsname}%
\fi
\fi\endgroup
\begin{picture}(1216,1230)(8693,-3368)
\end{picture}

%% file: nastyface.tex
\chapter{The Hidden Faces}\label{sec:nastyface}

In Chapters~\ref{sec:vassiliev} and~\ref{sec:tying}, I discussed in
some detail how to form an integral with trivalent graphs in which
the principal faces cancel.  In this chapter I'll show that, for
almost all the remaining faces of the fibered chain for a trivalent
graph $\Gamma$, the restriction to that face is 0.  The proofs are
mostly after Bott and Taubes~\cite{bt:knot}.

\section {Preliminaries}

A general principle that will be useful is the notion of
\emph{pullbacks} for bundles.  Given three bundles $B_1, B_2, B_3$
over $\knots$ and fibration preserving maps $i_1: B_1 \rightarrow B_3$
and $i_2: B_2 \rightarrow B_3$, we can ``complete the square'' in a
anonical way. We construct a bundle $B$ over $\knots$ with fiber
maps $j_1: B
\rightarrow B_1$ and $j_2: B \rightarrow B_2$ so that the following
square commutes:
\begin{equation}\begin{CD}
B @>>{j_1}> B_1\\
@VV{j_2}V     @VV{i_1}V\\
B_2 @>>{i_2}> B_3
\end{CD}\end{equation}
and $B$ is minimal: given any other $B'$ and maps
$j_1'$, $j_2'$, there is a unique map $\kappa: B' \rightarrow B$
so that $j_1 \circ \kappa = j_1'$ and $j_2 \circ \kappa = j_2'$.
Concretely, $B$ is the bundle
\begin{equation}
B = \{(x \in B_1, y \in B_2): i_1(x) = i_2(y)\}
\end{equation}

For example, for any finite set $A$ and $B, C \subset A$ with $B \cup
C = A$, there is a
commutative square
\begin{equation}\begin{CD}
C_A(\RR^3) @>>> C_B(\RR^3)\\
@VVV            @VVV\\
C_C(\RR^3) @>>> C_{B \cap C}(\RR^3)
\end{CD}\end{equation}
and $C_A(\RR^3)$ is a pullback in this diagram.  (Each of the maps is
a version of the natural project $C_A(\RR^3) \rightarrow C_B(\RR^3)$
for $B \subset A$, forgetting the points in $A \setminus B$.)
The spaces $C_{A;B}$ can be constructed in this way: they are the
pullbacks in
\begin{equation}\begin{CD}
C_{A;B} @>>> C_B(S^1)\\
@VVV          @VV{C_B(K)}V\\
C_A(\RR^3) @>>> C_B(\RR^3)
\end{CD}\end{equation}
This is how $C_{A;B}$ was defined in~\cite{bt:knot}, but it's a little bit
tricky to show that it is a manifold with corners in this way.

As an application, we have the following proposition.

\begin{proposition} \label{prop:notri}
The fibered chain corresponding to any trivalent graph
$\Gamma$ with
a subset $A \subset V_i(\Gamma)$ of the internal vertices of $\Gamma$
such that
\begin{itemize}
\item $|A| \ge 2$ and
\item The number of edges of $\Gamma$ connecting a vertex in $A$ to a
vertex not in $A$ is $\le 3$.
\end{itemize}
is degenerate.
\end{proposition}
\begin{proof}
Let $B$ be the set of vertices not in $A$ that are joined to a vertex
in $A$.  By hypothesis, $|B| \le 3$.
Now, $C_{V;c}$ (the configuration space appropriate to this
graph) is the pullback in
\begin{equation}\begin{CD}
C_{V;c} @>>> C_{A\cup B}(\RR^3)\\
@VVV @VVV\\
C_{(V\setminus A);c} @>>> C_B(\RR^3)
\end{CD}\end{equation}
and each summand in the fibered chain $\fchain(\Gamma)$ will be the
restriction of a product of two maps, from $C_{A\cup B}(\RR^3)$ and
$C_{(V\setminus A);c}$, to $C_{V;c} \subset C_{A\cup B}(\RR^3) \times
C_{(V\setminus A);c}$.  The two maps corespond, respectively, to the
graph $\Gamma_A$ whose edges are all edges in $\Gamma$ meeting $A$ and
the graph $\Gamma_{\bar A}$
containing all other edges of $\Gamma$.

But the first map (corresponding to $\Gamma_A$) lands in a proper
submanifold of $(S^2)^{E(\Gamma_A)}$.  The source, $C_{A\cup
B}(\RR^3)$, has dimension $3(|A| + |B|)$, while the
target, $(S^2)^{E(\Gamma_A)}$, has dimension
$2|E(\Gamma_A)| \ge 3|A|+|B|+d$ ($d$ is the number of vertices of
$B$ connected to more than one vertex of $A$), so this isn't
immediately obvious.
But the directions between all these points is unchanged under
\begin{itemize}
\item overall translation of $A$ and $B$ in $\RR^3$ (3 dimensions)
\item moving a vertex in $B$ away from the vertex it's connected to in
$A$, for vertices in $B$ connected to exactly one vertex of $A$ ($|B|
- d$ dimensions)
\item overall scaling (1 dimension)
\end{itemize}
All of these transformations are independent because $|A| \ge 2$.
Taking the quotient by all these transformations, we see that the map
from $C_{A\cup B}(\RR^3)$ factors through a space of dimension
$3(|A|+|B|)-3-|B|+d-1 \le 3|A|+|B|+d-1 < \dim (S^2)^{E(\Gamma_A)}$, so
indeed we do land in a proper submanifold.

The homology chain $\fchain(\Gamma)$ will be inside the product of
this submanifold with $(S^2)^{E(\Gamma_{\bar A})}$, and thus in a
proper submanifold of
$(S^2)^{E(\Gamma)}$.
\end{proof}

Arguments along similar lines show that boundary components
$C_{(A;B),D}$ are
degenerate if they have enough (or too few) \emph{external edges}:
edges connecting a point in $D$ with a point not in $D$.  (For
trivalent diagrams, boundary components with $B\cap D =\emptyset$ are
degenerate if there are more than 4 external edges, while boundary
components with $B\cap D \ne \emptyset$ are degenerate if there are
more than 2 external edges.)  See~\cite{bt:knot} for details.
See below for a proof 
by global symmetry (or sometimes degeneracy) that these faces vanish
as chains.

\section{Conditions for Vanishing}

There are three different arguments to show that a boundary face vanishes.
Throughout, we're considering the face of the configuration space
$C_{V,c}$ corresponding to the trivalent graph $\Gamma$ in which
points in the subset $A$ collide.  Since these are non-principal
faces, by definition $|A| \ge 3$.
$\Gamma_A$ is the subgraph of
$\Gamma$ whose edges all propogators that connect points of $A$ (note
that this differs from the definition in the proof of
Proposition~\ref{prop:notri}), and
$\Gamma_{\bar A}$ consists of all other edges of $\Gamma$.  $M =
(S^2)^{E(\Gamma)}$ is the target for the fibered homology chain
$\fchain(\Gamma)$.

\begin{case}\label{case:notconn}
If $\Gamma_A$ is not connected, the fibered chain
$\fchain(C_{(V;c),A})$ is degenerate. 
\end{case}
\begin{proof}
In this case, we can independently translate the vertices of two
disjoint subsets of $\Gamma_A$ within
$C_{A;(c\cap A)}(\calT\,\RR^3;\dot K(s_A))$ 
without changing the map to $M$.  The
translation is in any direction in $\RR^3$ (if $A \cap c = \emptyset$
or along the direction $\dot K(s_A)$ (otherwise), but in either case
the map to $M$ factors through a map to a bundle of strictly lower
dimension.  The quotient is non-trivial since $|A| \ge 3$.  (Recall that
the local coordinates are already the quotient by translation and
scaling of all the points in $A$.)

Note that this same argument works even if $|A| = 2$ but both vertices
are free.  This justifies ignoring certain principal faces.
\end{proof}

\begin{case}\label{case:oneext}
If $\Gamma$ has a subgraph like
\begin{equation}
\mathcenter{\input{draws/subgraph-1.tex}}
\end{equation}
(with $a,b,c \in A$ and $d \notin A$) then the fibered chain
$\fchain(C_{(V;c),A})$ 
vanishes by a global symmetry.
\end{case}
\begin{proof}
The map on
$C_{A;(c\cap A)}(\calT\,\RR^3;\dot K(s_A))$ 
that takes $x_a$ to $x_b + x_c - x_a$ and leaves all other coordinates
unchanged is an orientation-reversing
automorphism on the stratum%
\footnote{The pedantic reader will note that this is not quite true,
since $x_b + x_c - x_a$ might happen to
coincide with $x_e$ for some other point $e \in A$.  For this reader,
I will point out that we didn't actually need to blow up on the
diagonal $\Delta_{ae}$, since $a$ and $e$ cannot be connected by an
edge of $\Gamma$, and if we don't do this blowup we do get a true
automorphism.}
that takes the fibered chain to itself (swapping and negating the spheres
corresponding to the edges $a$--$b$ and $a$--$c$).
\end{proof}

\begin{case}\label{case:twoext}
If $\Gamma$ has a subgraph like
\begin{equation}
\mathcenter{\input{draws/subgraph-2.tex}}
\end{equation}
(with $a,b,c,d \in A$ and $e,f \notin A$) then the fibered chain
$\fchain(C_{(V;c),A})$ 
vanishes by a global symmetry.
\end{case}
\begin{proof}
The map on
$C_{A;(c\cap A)}(\calT\,\RR^3;\dot K(s_A))$ 
that takes $x_b$ to $x_c + x_d - x_b$, $x_a$ to $x_c + x_d - x_a$, and
leaves all other coordinates fixed is an orientation-preserving
autmorphism on the stratum\footnote{See previous footnote} that takes
the fibered chain to its negative (swapping and negating the spheres
corresponding to the edges $b$--$c$ and $b$--$d$ and negating the
sphere corresponding to $a$--$b$).
\end{proof}

\section{Most Hidden Faces Vanish}
\begin{proposition}
Any non-principal face with external edges vanishes as a fibered chain.
\end{proposition}
\begin{proof}
Suppose the face corresponding to the subset $A$ has external edges.
If some vertex $a\in A$ is connected to exactly one external vertex,
either $a$ is on the knot, in which case $\Gamma_A$ (as in the
previous section) cannot be connected and the face vanishes by
Case~\ref{case:notconn}; or $a$ is not on the not, in which case
$\Gamma$ has a subgraph as in Case~\ref{case:oneext}.  If some vertex
$a$ is connected to three external vertices, $\Gamma_A$ cannot be
connected and the face is degenerate by Case~\ref{case:notconn}.  So
suppose boundary vertices (vertices meeting at least one external
edge) are meet exactly two external edges.  If some boundary vertex
$a$ is also connected to a vertex $b \in A$ on the knot or another
boundary vertex, then (since
there is at least one other vertex in $A$), $\Gamma_A$ is not
connected and the face is degenerate.  Otherwise, $a$ must be
connected to another non-boundary vertex not on the knot and $\Gamma$
will have a subgraph as in Case~\ref{case:twoext}.
\end{proof}

There are two cases left: faces in which some points $A \subset
V(\Gamma)$ approach $\infty$, and faces in which some points $A$,
with no external edges, approach each other.  The first case is easy
(by Proposition~\ref{prop:notri}, if $\fchain(\Gamma)$ is
non-degnerate $A$ will have at least three
external edges, all of which map to the same point in $S^2$, so the
image of the map has codimension at least 4; since the face has
codimension 1, the map is degenerate).  The second is not so easy, and
in fact, it's not always true that these faces vanish.  They are
called \emph{anomalous faces}, and are the subject of the next
chapter.


%% file: draws/subgraph-1.tex
\begin{picture}(0,0)%
\epsfig{file=draws/subgraph-1.pstex}%
\end{picture}%
\setlength{\unitlength}{3947sp}%
\begingroup\makeatletter\ifx\SetFigFont\undefined
\def\x#1#2#3#4#5#6#7\relax{\def\x{#1#2#3#4#5#6}}%
\expandafter\x\fmtname xxxxxx\relax \def\y{splain}%
\ifx\x\y   
\gdef\SetFigFont#1#2#3{%
  \ifnum #1<17\tiny\else \ifnum #1<20\small\else
  \ifnum #1<24\normalsize\else \ifnum #1<29\large\else
  \ifnum #1<34\Large\else \ifnum #1<41\LARGE\else
     \huge\fi\fi\fi\fi\fi\fi
  \csname #3\endcsname}%
\else
\gdef\SetFigFont#1#2#3{\begingroup
  \count@#1\relax \ifnum 25<\count@\count@25\fi
  \def\x{\endgroup\@setsize\SetFigFont{#2pt}}%
  \expandafter\x
    \csname \romannumeral\the\count@ pt\expandafter\endcsname
    \csname @\romannumeral\the\count@ pt\endcsname
  \csname #3\endcsname}%
\fi
\fi\endgroup
\begin{picture}(2475,1383)(10147,-5674)
\put(11043,-5141){\makebox(0,0)[lb]{\smash{\SetFigFont{12}{14.4}{rm}$a$}}}
\put(10424,-5674){\makebox(0,0)[rb]{\smash{\SetFigFont{12}{14.4}{rm}$b$}}}
\put(12005,-5181){\makebox(0,0)[lb]{\smash{\SetFigFont{12}{14.4}{rm}$A$}}}
\put(12005,-4716){\makebox(0,0)[lb]{\smash{\SetFigFont{12}{14.4}{rm}not $A$}}}
\put(10948,-4456){\makebox(0,0)[b]{\smash{\SetFigFont{12}{14.4}{rm}$d$}}}
\put(11500,-5674){\makebox(0,0)[lb]{\smash{\SetFigFont{12}{14.4}{rm}$c$}}}
\end{picture}

%% file: draws/subgraph-2.tex
\begin{picture}(0,0)%
\epsfig{file=draws/subgraph-2.pstex}%
\end{picture}%
\setlength{\unitlength}{3947sp}%
\begingroup\makeatletter\ifx\SetFigFont\undefined
\def\x#1#2#3#4#5#6#7\relax{\def\x{#1#2#3#4#5#6}}%
\expandafter\x\fmtname xxxxxx\relax \def\y{splain}%
\ifx\x\y   
\gdef\SetFigFont#1#2#3{%
  \ifnum #1<17\tiny\else \ifnum #1<20\small\else
  \ifnum #1<24\normalsize\else \ifnum #1<29\large\else
  \ifnum #1<34\Large\else \ifnum #1<41\LARGE\else
     \huge\fi\fi\fi\fi\fi\fi
  \csname #3\endcsname}%
\else
\gdef\SetFigFont#1#2#3{\begingroup
  \count@#1\relax \ifnum 25<\count@\count@25\fi
  \def\x{\endgroup\@setsize\SetFigFont{#2pt}}%
  \expandafter\x
    \csname \romannumeral\the\count@ pt\expandafter\endcsname
    \csname @\romannumeral\the\count@ pt\endcsname
  \csname #3\endcsname}%
\fi
\fi\endgroup
\begin{picture}(2454,1725)(10162,-6046)
\put(12005,-5181){\makebox(0,0)[lb]{\smash{\SetFigFont{12}{14.4}{rm}$A$}}}
\put(12005,-4716){\makebox(0,0)[lb]{\smash{\SetFigFont{12}{14.4}{rm}not $A$}}}
\put(10351,-4486){\makebox(0,0)[b]{\smash{\SetFigFont{12}{14.4}{rm}$e$}}}
\put(11551,-4486){\makebox(0,0)[b]{\smash{\SetFigFont{12}{14.4}{rm}$f$}}}
\put(11101,-5611){\makebox(0,0)[lb]{\smash{\SetFigFont{12}{14.4}{rm}$b$}}}
\put(10424,-6031){\makebox(0,0)[rb]{\smash{\SetFigFont{12}{14.4}{rm}$c$}}}
\put(11470,-6046){\makebox(0,0)[lb]{\smash{\SetFigFont{12}{14.4}{rm}$d$}}}
\put(11101,-5236){\makebox(0,0)[lb]{\smash{\SetFigFont{12}{14.4}{rm}$a$}}}
\end{picture}

%% file: anomaly.tex
\chapter{Anomalous Faces}
\label{sec:anomaly}
In the previous section, I showed that most hidden faces (faces that
aren't just two vertices coming together) vanish.  The remaining
faces, dubbed ``anomalous faces'' by Bott and Taubes in~\cite{bt:knot},
are much more interesting.  They do \emph{not} all vanish as homology 
chains.  Fortunately, most of those anomalous faces that do not vanish can be
corrected by an additional term.




\section{Preliminary Remarks}

First, some definitions and reductions of the problem.

\begin{definition}
A \emph{split} graph $\Gamma$ is a graph whose vertices can be
partitioned into two sets, $A$ and $B$, such that the only edges
running between $A$ and $B$ are cycle (knot) edges.
\end{definition}

The notion of a split graph is closely related to the notion of
connected sum of knots.

\begin{definition}
The \emph{connected sum} of two knots $\A$ and $\B$ is the knot formed
by embedding $\A$ and $\B$ in $\RR^3$ so that $\A$ is contained in a
ball that does not intersect $\B$,
removing a short section of $\A$ and a short section of $\B$ and
connecting them with a strand from $\A$ to $\B$ and a strand from $\B$
to $\A$ so that the orientations match and the loop formed by the two
sections removed and the two strands added is an unknot unlinked with $\A$ and
with $\B$.  See Figure~\ref{fig:connectsum} for an example.  For a
proof that the isotopy class of the resulting knot does not depend on
the choices, see~\cite{kauffman:knotsphysics}.
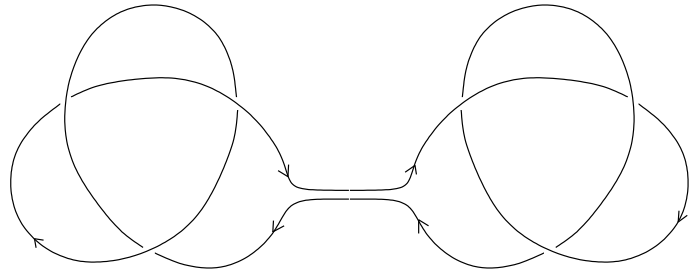
\begin{figure}
\centerline{\input{draws/connect-sum.tex}}
\caption{Connected sum of two knots}
\label{fig:connectsum}
\end{figure}
\end{definition}

A split trivalent graph can be thought of as the connected sum of two
trivalent graphs
(although this operation depends on the base point chosen; but see
Proposition~\ref{prop:connsumweight} below).
A weight system with no split graphs (called a \emph{primitive} weight
system) corresponds to an invariant (to
be constructed)
whose value on a connected sum of two knots is the sum of the values
on the individual knots.  For the corresponding homology semicycle,
the same thing is true, at least if the two knots are widely separated
and the strands connecting them are near each other.

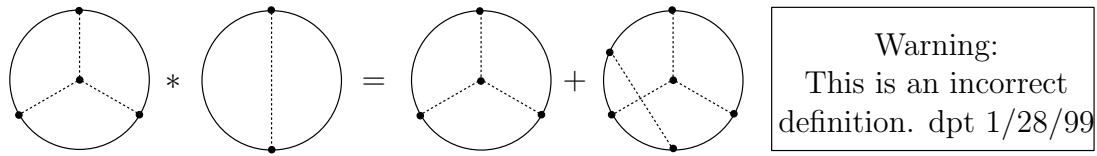
\begin{figure}
\centerline{\input{draws/graph-prod.tex}}
\caption{The product of graphs}
\label{fig:graphprod}
\end{figure}

There is a natural product on the vector space spanned by trivalent
graphs: the sum of all ways of
interleaving the knot vertices of the two graphs, as in
Figure~\ref{fig:graphprod}.  (The orientations carry over as well.)
This product clearly corresponds to taking the product of the
corresponding counting formulas or integrals.  The product of two
weight systems will again be a weight system.  (The STU relation
between two knot vertices from different original graphs is the only
new 
thing to check, and it's satisfied because of the sum over all
interleavings.)

Note that the product of two non-trivial weight systems will not be
primitive.  Conversely, we have the following proposition, which
corresponds nicely to the (unique) decomposition of a knot into pieces
that cannot be expressed as a connected sum (``prime'' knots).

\begin{proposition}
\label{prop:connsumweight}
Every weight system can be expressed uniquely as a sum of products of
primitive weight systems.
\end{proposition}
\begin{proof}
See~\cite{bn:vassiliev}, and especially Proposition 3.4.
\end{proof}

Since it's obviously desireable to have the product of weight systems
correspond to the product of the (numerical) invariants, it suffices
to construct invariants for only the primitive weight diagrams.  For
these diagrams, the anomalous faces must have the entire graph
collapsing to a single point on the knot.  It's interesting to note
that there are several connected components to this face for a given graph,
corresponding to different ways to make the cyclic order of vertices
into a linear order.

One note: the anomalous face is degenerate if the graph consisting of just
the propogators is not connected, by Case~\ref{case:notconn} of the
previous section.  This means that any graph with a non-trivial
anomalous face has a relatively large number of internal vertices.

For a given non-split, connected graph $\Gamma$, the restriction of
$\fchain(\Gamma)$ to the (unique) anomalous stratum defines a map from
the appropriate configuration space, $C_{(V;c),V}$, to
$M = (S^2)^{E(\Gamma)}$.  Since all the points collapse, the
configuration space is fibered over $C_1(K) \simeq S^1$.  Each
component of the fiber is $C_{V;c}(\calT\,\RR^3;\dot K(s_V))$ (where
$c(\Gamma)$ is given a linear order compatible with its cyclic order).
This can be viewed as a bundle $B$ over $S^2$, with projection map
$\pi: B \rightarrow S^2$; then
$C_{(V;c),V}$ is the pullback to $C_1(K)$ of $B_{V;c}$ via the
(normalized) tangent
map.  This is a different kind of pullback than we saw before.  This
kind of pullback can be constructed whenever you have a bundle $B$
over a manifold $Y$ (with projection map $\pi: B \rightarrow Y$) and a
map $f: X \rightarrow Y$ from another
manifold $X$; then the pullback of the bundle $B$ to $X$ is the
universal $A$ in
\begin{equation}\begin{CD}
A @>>> B\\
@VVV @VV{\pi}V\\
X @>>{f}> Y
\end{CD}\end{equation}
or, concretely, the manifold $A = \{(b \in B, x \in X): f(x) =
\pi(b)\}$.

The sum of maps $\fchain(\Gamma): C_{V;c}(\calT\,\RR^3;\dot K(s_V))
\rightarrow M$ is
the pullback of a sum of maps $\phi: B \rightarrow M$.  $\phi$ gives a
codimension 2 fibered chain on $S^2$ that captures many of the
properties of $\fchain(\Gamma)$.
For instance, if $B(\pi^{-1}(\xi)) = 0$ as a chain (modulo singular
chains) in $M$ for each $\xi\in S^2$,
then $\fchain(\Gamma)$ restricted to the anomalous face is, like
$\phi$, a degenerate fibered chain.  Or in the integral formulation:
if the integral over the fibers of $\phi$ of the pullback of a volume
form on $M$ is the 2-form $\omega$ on $S^2$, the corresponding integral
over the fibers of $\fchain(\Gamma)$ is the integral of the pullback
of $\omega$ via the tangent map.  To show that the anomalous face
vanishes as a fibered chain, it suffices to show that $\phi$ does (or
equivalently, that the integral over the fibers of $\phi$ of an
$n$-form $\omega$ on $M$ is 0 for every $\omega$). 
To show that the integral constructed with the canonical form on $M$
is an invariant, it suffices to show that the integral over the fibers
of the canonical form $\theta$ on $M$ is the zero 2-form.

\section{Parity of Invariants}

There are two orientations involved in the definition of a knot: the
orientation on $S^1$ and the orientation on $\RR^3$.  Reversing the
orientation on $S^1$ yields the \emph{inverse} of the knot; reversing
the orientation on $\RR^3$ yields its \emph{reverse??}.  A (numerical)
knot 
invariant that takes the same value on a knot and its inverse is
called \emph{even} under reversal of $S^1$; invariants that take the
negative value on the inverse are called \emph{odd} under reversal of
$S^1$.  Even and odd invariants under reversal of $\RR^3$ are defined
similarly.  Every invariant can be uniquely expressed as a sum of
invariants 
that are either even or odd under reversal of $S^1$ and even or odd
under reversal of $\RR^3$.

Knots that
are not equal to their reverse are easy to come by.  The first
non-trivial knot, the trefoil, is one such.  Invariants odd under
reversal of $\RR^3$ are correspondingly easy to find.  Every weight
diagram of odd degree yields an invariant odd under reversal of
$\RR^3$.  We expect this to be true since the definition of derivative
depends on the orientation of $\RR^3$, so the derivative of an
invariant (its evaluation on an immersed loop with one more crossing)
changes its $\RR^3$ parity, and a constant (an invariant whose derivative is
0) is necessarily even under reversal of $\RR^3$.  And in fact, by
looking at the tinkertoy counting formula and
noting that each $\epsilon$-tensor in the sign algorithm changes sign
when you reverse the orientation of $\RR^3$, one can show that the
function on loops associated to a
trivalent graph changes sign under reversal of $\RR^3$ exactly when
the graph is
of odd degree (i.e., has $2k$ vertices, for $k$ odd).  Our correction
term below will not change this.

Knots that are not equal to their inverse
are much rarer.  The knot with the fewest number of crossings that is
not equal to its inverse has eight crossings.  Vassiliev invariants
that are odd under reversal of $S^1$ are also rare, and may not
exist at all.  In~\cite{bn:vassiliev}, Bar-Natan presents the results of
computer calculations enumerating the weight diagrams up to
degree 9; so far, none that are odd under reversal of $S^1$ have been
found.  (The definition of derivative of a knot invariant does not
depend on the orientation of the knot (check this!), so the $S^1$
parity of a knot invariant can be determined by looking at its largest
non-zero derivative, or a weight diagram, which is what Bar-Natan
computed.)

The parities of a (pure in parity, primitive) weight diagram $D$ give
important information about the
bundle map $\phi$ discussed above.  In particular, they tell us the
parity of $\phi$:
whether $\phi \circ N$ (with $N: S^2 \rightarrow S^2$ the antipodal
map) is $\pm \phi$.  If $D$ is of degree $k$, then $\phi \circ N =
(-1)^{k+1}\phi$, while if $D$ is of $S^2$ parity $\epsilon = \pm 1$,
then $\phi \circ N = \epsilon\phi$.  The easiest way to see this is to
note that the variation of a loop function of a given parity will have
the same parity (for either type of parity).  The integral of a 2-form
of a given $\RR^3$ parity over the knot will have the opposite $\RR^3$
parity (as a 1-form on $\knots$), since this integral depends on the
orientation on $S^2$, which is inherited from the orientation of
$\RR^3$.  The integral of a 2-form of a given $S^1$ parity over the
knot has the same parity: the integral depends on the orientation on
$S^1$, but it also depends on the tangent map $\dot K$, which is
reversed (an orientation-reversing map on the sphere) when the
orientation on $S^1$ is changed.  (It's possible to exhibit
explicit automorphisms of the bundle $B$ showing both these facts.)

\begin{theorem}
There exists a canonically-defined fibered homology cycle for each
weight diagram $D$, with a correction term if $D$ has odd $\RR^3$
parity and even $S^1$ parity.  In any case, no correction term is
necessary for almost-planar knots or for the canonical integral.
\end{theorem}
\begin{proof}
It's natural
to proceed by cases, depending on the parities of the weight diagram
$D$.

\setcounter{case}{1}
\begin{case}
$\RR^3$ parity even, $S^1$ parity even or $\RR^3$ parity odd, $S^1$
parity odd
\end{case}
\begin{proof}
In either of these two cases, the fibered chain $\phi$ must be both
odd and even under the antipode map on $S^2$ and hence vanishes.
\end{proof}

\begin{case}
$\RR^3$ parity odd, $S^1$ parity even
\end{case}
\begin{proof}
This is the only case in which we'll need a correction factor.  The
correction will be given by $-\frac{1}{2}$ the pullback of the bundle
$B$ and map
$\phi$ to the space $C_2(K)$ via the map $\frac{x_2 - x_1}{|x_2-x_1|^3}$.
$C_2(K)$ has a single boundary stratum with two connected components
(the points 1 and 2 coming together in either of the two possible
orders).  The variations of these two faces
are the pullbacks to $C_1(K)$ via the tangent map of
$-\frac{1}{2} \phi$ and $-\frac{1}{2} \phi \circ N$.  Since
$\phi\circ N = \phi$, these two terms are the same and combine to
cancel the original variation, yielding a fibered homology cycle and
hence a knot invariant.

The tinkertoy structures for this correction term are a little more
involved than usual.  They consist of two nodes on the knot and a
``virtual'' rod
between them (a rod \emph{not} constrained to lie in one of the chosen
directions), together with a tinkertoy structure (for the graph whose
anomaly we're correcting) constructed \emph{on the virtual rod}; i.e.,
each node that would have gone on the knot in a regular tinkertoy
diagram must lie on the virtual rod.  For general sets of chosen
directions, this condition can only be satisfied if the virtual rod is
in a finite number of positions.
Given any such ``virtual''
tinkertoy diagram, its clear that translations (along the rod) and
scalings of it will also satisfy the conditions; these should only be
counted once.  The sign will also be more involved.  The only property
of the sign I'll need is that (like all bundles over $C_2(K)$ in this
way) the sign of an overcrossing and an undercrossing (with the same
virtual diagram, relative to the corresponding virtual direction) are
opposite.

This correction term takes a special form for almost planar knots.  An
\emph{almost planar knot} is, as the name implies, very close to lying
in a plane (over which the knot should have generic crossings).
Although only the unknot can exactly lie in a plane, any
knot can come arbitrarily close to a plane, with the only deviations
from true planarity near the crossings of the plane projection.

For a knot in such a position, the virtual directions (and hence a
correction term) can only occur at the crossings, since, relative to
the almost planar knot, these directions are almost vertical (for
generic initially chosen directions).  Furthermore, at each crossing
every virtual direction occurs (some with point 1 on the lower strand
and point 2 on the upper strand,
some vice versa).  The net contribution from each overcrossing will
therefore be the same, say $\lambda$ and the net contribution from
each undercrossing will be $-\lambda$.  The correction term for almost
planar
knots is then $\lambda$(\#overcrossings - \#undercrossings).  This
difference between the number of overcrossings and undercrossings of
an almost planar knot is called the \emph{writhe}; it is the framing
number of the framing given by choosing normal vectors parallel to the
plane.  (See~\cite{kauffman:knotsphysics}).

One can show $\lambda$ does not depend on the form $\omega$ chosen on
the product of $S^2$: this boils down to the familiar fact that the
boundary of a boundary is 0.

From field theory arguments one expects this correction term $\lambda$
to always be 0, but I don't know any proof.
The only non-trivial case which I know is that for the unique degree 3
Vassiliev invariant (shown in
Figure~\ref{fig:deg3})
\begin{figure}
\begin{equation}
\begin{split}
w\left(\mathcenter{\input{draws/deg3-1.tex}}\right) &= 2\qquad
  w\left(\mathcenter{\input{draws/deg3-2.tex}}\right) = -1\\
w\left(\mathcenter{\input{draws/deg3-3.tex}}\right) &= 1\qquad
  w\left(\mathcenter{\input{draws/deg3-4.tex}}\right) = 1
\end{split}
\end{equation}
\caption{The unique degree 3 weight system}
\label{fig:deg3}
\end{figure}
the only anomaly term comes from the last factor and I have found a
direct argument that $\lambda$ is 0 for this diagram, although the
correction term does not vanish for all forms $\omega$.

Since the canonical form on $(S^2)^e$ is SO(3) invariant, the
correction form will also be SO(3) invariant, and hence will vanish if
$\lambda$, its integral over the sphere, is 0.
\end{proof}

\begin{case}
$\RR^3$ parity even, $S^1$ parity odd
\end{case}
\begin{proof}
In this case, the variation form $\psi$ is odd under the antipodal
map $N$, and thus has total integral 0 over the sphere.  This implies
that $\psi$ is cohomologous to 0, or equivalently $\psi = d\alpha$ for
some 1-form $\alpha$ on $S^2$.

\begin{proposition}
If $d\alpha = \psi$, then $I(\Gamma)(K) - \int_{S^1} K^* \alpha$ is a
knot invariant that does not depend form $\alpha$ chosen, and is
independent of 
\end{proposition}

Now consider the integral $\int_{S^1} K^*\alpha$.  By
Proposition~\ref{prop:dpushforward}, the variation of this integral is
$\int_{S^1} K^*d\alpha = \int_{S^1} K^* \psi$ is exactly the variation
of our original form, so subtracting this term yields a knot
invariant.

although not cancellation formula of original type, doesn't change as
vary the knot.

For plane curves, multiple of the winding number.  (0 for canonical
integral.)

Surprising if it vanishes for all $\omega$, because then it
vanishes identically; but we've already looked at the obvious
symmetries that would keep direction fixed.

\end{proof}
This completes the proof of the proposition, and hence our main theorem.
\end{proof}

%% file: draws/connect-sum.tex
\begin{picture}(0,0)%
\epsfig{file=draws/connect-sum.pstex}%
\end{picture}%
\setlength{\unitlength}{3947sp}%
\begingroup\makeatletter\ifx\SetFigFont\undefined
\def\x#1#2#3#4#5#6#7\relax{\def\x{#1#2#3#4#5#6}}%
\expandafter\x\fmtname xxxxxx\relax \def\y{splain}%
\ifx\x\y   
\gdef\SetFigFont#1#2#3{%
  \ifnum #1<17\tiny\else \ifnum #1<20\small\else
  \ifnum #1<24\normalsize\else \ifnum #1<29\large\else
  \ifnum #1<34\Large\else \ifnum #1<41\LARGE\else
     \huge\fi\fi\fi\fi\fi\fi
  \csname #3\endcsname}%
\else
\gdef\SetFigFont#1#2#3{\begingroup
  \count@#1\relax \ifnum 25<\count@\count@25\fi
  \def\x{\endgroup\@setsize\SetFigFont{#2pt}}%
  \expandafter\x
    \csname \romannumeral\the\count@ pt\expandafter\endcsname
    \csname @\romannumeral\the\count@ pt\endcsname
  \csname #3\endcsname}%
\fi
\fi\endgroup
\begin{picture}(4278,1677)(9188,-4186)
\end{picture}

%% file: draws/graph-prod.tex
\begin{picture}(0,0)%
\epsfig{file=draws/graph-prod.pstex}%
\end{picture}%
\setlength{\unitlength}{3947sp}%
\begingroup\makeatletter\ifx\SetFigFont\undefined
\def\x#1#2#3#4#5#6#7\relax{\def\x{#1#2#3#4#5#6}}%
\expandafter\x\fmtname xxxxxx\relax \def\y{splain}%
\ifx\x\y   
\gdef\SetFigFont#1#2#3{%
  \ifnum #1<17\tiny\else \ifnum #1<20\small\else
  \ifnum #1<24\normalsize\else \ifnum #1<29\large\else
  \ifnum #1<34\Large\else \ifnum #1<41\LARGE\else
     \huge\fi\fi\fi\fi\fi\fi
  \csname #3\endcsname}%
\else
\gdef\SetFigFont#1#2#3{\begingroup
  \count@#1\relax \ifnum 25<\count@\count@25\fi
  \def\x{\endgroup\@setsize\SetFigFont{#2pt}}%
  \expandafter\x
    \csname \romannumeral\the\count@ pt\expandafter\endcsname
    \csname @\romannumeral\the\count@ pt\endcsname
  \csname #3\endcsname}%
\fi
\fi\endgroup
\begin{picture}(6830,934)(9833,-4327)
\put(15676,-3705){\makebox(0,0)[b]{\smash{\SetFigFont{12}{14.4}{rm}Warning:}}}
\put(15676,-3945){\makebox(0,0)[b]{\smash{\SetFigFont{12}{14.4}{rm}This is an incorrect}}}
\put(15676,-4185){\makebox(0,0)[b]{\smash{\SetFigFont{12}{14.4}{rm}definition. dpt 1/28/99}}}
\put(10882,-3915){\makebox(0,0)[b]{\smash{\SetFigFont{12}{14.4}{rm}$*$}}}
\put(12143,-3915){\makebox(0,0)[b]{\smash{\SetFigFont{12}{14.4}{rm}$=$}}}
\put(13404,-3915){\makebox(0,0)[b]{\smash{\SetFigFont{12}{14.4}{rm}$+$}}}
\end{picture}

%% file: draws/deg3-1.tex
\begin{picture}(0,0)%
\epsfig{file=draws/deg3-1.pstex}%
\end{picture}%
\setlength{\unitlength}{3947sp}%
\begingroup\makeatletter\ifx\SetFigFont\undefined
\def\x#1#2#3#4#5#6#7\relax{\def\x{#1#2#3#4#5#6}}%
\expandafter\x\fmtname xxxxxx\relax \def\y{splain}%
\ifx\x\y   
\gdef\SetFigFont#1#2#3{%
  \ifnum #1<17\tiny\else \ifnum #1<20\small\else
  \ifnum #1<24\normalsize\else \ifnum #1<29\large\else
  \ifnum #1<34\Large\else \ifnum #1<41\LARGE\else
     \huge\fi\fi\fi\fi\fi\fi
  \csname #3\endcsname}%
\else
\gdef\SetFigFont#1#2#3{\begingroup
  \count@#1\relax \ifnum 25<\count@\count@25\fi
  \def\x{\endgroup\@setsize\SetFigFont{#2pt}}%
  \expandafter\x
    \csname \romannumeral\the\count@ pt\expandafter\endcsname
    \csname @\romannumeral\the\count@ pt\endcsname
  \csname #3\endcsname}%
\fi
\fi\endgroup
\begin{picture}(932,894)(9436,-4408)
\end{picture}

%% file: draws/deg3-2.tex
\begin{picture}(0,0)%
\epsfig{file=draws/deg3-2.pstex}%
\end{picture}%
\setlength{\unitlength}{3947sp}%
\begingroup\makeatletter\ifx\SetFigFont\undefined
\def\x#1#2#3#4#5#6#7\relax{\def\x{#1#2#3#4#5#6}}%
\expandafter\x\fmtname xxxxxx\relax \def\y{splain}%
\ifx\x\y   
\gdef\SetFigFont#1#2#3{%
  \ifnum #1<17\tiny\else \ifnum #1<20\small\else
  \ifnum #1<24\normalsize\else \ifnum #1<29\large\else
  \ifnum #1<34\Large\else \ifnum #1<41\LARGE\else
     \huge\fi\fi\fi\fi\fi\fi
  \csname #3\endcsname}%
\else
\gdef\SetFigFont#1#2#3{\begingroup
  \count@#1\relax \ifnum 25<\count@\count@25\fi
  \def\x{\endgroup\@setsize\SetFigFont{#2pt}}%
  \expandafter\x
    \csname \romannumeral\the\count@ pt\expandafter\endcsname
    \csname @\romannumeral\the\count@ pt\endcsname
  \csname #3\endcsname}%
\fi
\fi\endgroup
\begin{picture}(899,918)(9453,-4428)
\end{picture}

%% file: draws/deg3-3.tex
\begin{picture}(0,0)%
\epsfig{file=draws/deg3-3.pstex}%
\end{picture}%
\setlength{\unitlength}{3947sp}%
\begingroup\makeatletter\ifx\SetFigFont\undefined
\def\x#1#2#3#4#5#6#7\relax{\def\x{#1#2#3#4#5#6}}%
\expandafter\x\fmtname xxxxxx\relax \def\y{splain}%
\ifx\x\y   
\gdef\SetFigFont#1#2#3{%
  \ifnum #1<17\tiny\else \ifnum #1<20\small\else
  \ifnum #1<24\normalsize\else \ifnum #1<29\large\else
  \ifnum #1<34\Large\else \ifnum #1<41\LARGE\else
     \huge\fi\fi\fi\fi\fi\fi
  \csname #3\endcsname}%
\else
\gdef\SetFigFont#1#2#3{\begingroup
  \count@#1\relax \ifnum 25<\count@\count@25\fi
  \def\x{\endgroup\@setsize\SetFigFont{#2pt}}%
  \expandafter\x
    \csname \romannumeral\the\count@ pt\expandafter\endcsname
    \csname @\romannumeral\the\count@ pt\endcsname
  \csname #3\endcsname}%
\fi
\fi\endgroup
\begin{picture}(896,908)(9453,-4418)
\end{picture}

%% file: draws/deg3-4.tex
\begin{picture}(0,0)%
\epsfig{file=draws/deg3-4.pstex}%
\end{picture}%
\setlength{\unitlength}{3947sp}%
\begingroup\makeatletter\ifx\SetFigFont\undefined
\def\x#1#2#3#4#5#6#7\relax{\def\x{#1#2#3#4#5#6}}%
\expandafter\x\fmtname xxxxxx\relax \def\y{splain}%
\ifx\x\y   
\gdef\SetFigFont#1#2#3{%
  \ifnum #1<17\tiny\else \ifnum #1<20\small\else
  \ifnum #1<24\normalsize\else \ifnum #1<29\large\else
  \ifnum #1<34\Large\else \ifnum #1<41\LARGE\else
     \huge\fi\fi\fi\fi\fi\fi
  \csname #3\endcsname}%
\else
\gdef\SetFigFont#1#2#3{\begingroup
  \count@#1\relax \ifnum 25<\count@\count@25\fi
  \def\x{\endgroup\@setsize\SetFigFont{#2pt}}%
  \expandafter\x
    \csname \romannumeral\the\count@ pt\expandafter\endcsname
    \csname @\romannumeral\the\count@ pt\endcsname
  \csname #3\endcsname}%
\fi
\fi\endgroup
\begin{picture}(927,903)(9434,-4408)
\end{picture}

%% file: linalg.tex
\chapter{Some linear algebra} \label{sec:linalg}

In this section I will present some results in linear algebra, related
to the determinant, that may not be familiar to all readers.  I assume
a basic knowledge of linear algebra, including tensor products and
antisymmetric products and the exterior algebra.

The \emph{determinant} of an $n$\/-dimensional vector space $V$ is the
1-dimensional vector space
\begin{equation}
\det V = \bigwedge^{\text{top}} V = \bigwedge^n V
\end{equation}
This is a functor from the category $\Vect^{\text{i}}$
of vector spaces and isomorphisms
to itself.  ($\bigwedge^n$ is a functor on the ordinary category of
vector spaces and linear maps.  The restriction to isomorphisms
ensures that the domain and range have the same dimension.)
It's called the determinant because the ordinary determinant of a linear map
$f: V \rightarrow V$ is the (unique) eigenvalue of the linear
map $\bigwedge^n f: \bigwedge^n V \rightarrow \bigwedge^n V$.

\begin{lemma} \label{lem:basicisom}
$\det V \otimes \det W$ is canonically isomorphic to $\det (V
\oplus W)$.
\end{lemma}
\begin{proof} These are both 1-dimensional vector spaces, so a canonical
non-zero map
between them will suffice.  In general, a linear map from $A \otimes
B$ to $C$ is a bilinear map from $A$ and $B$ to $C$.  So suppose $a
\in \det V$ and $b \in \det W$.  Since $V$ and $W$ are both subspaces
of $V \oplus W$, $\det V$ is a subspace of $\bigwedge^{\dim V} (V
\oplus W)$ and $\det W$ is a subspace of $\bigwedge^{\dim W} (V
\oplus W)$.  If we consider both $a$ and $b$ in $\bigwedge^* (V \oplus
W)$, then $a \wedge b$ is the desired element of $\det (V \oplus W)$.
This is a bilinear map and is clearly functorial.
\end{proof}

Note that the isomorphism in the above lemma depends on the order $V$
and $W$ are specified.  Specifically, we have the following:

\begin{lemma} \label{lem:rearr}
Under the isomorphisms
\begin{equation}
\det V \otimes \det W \simeq \det(V\oplus W)\simeq\det(W\oplus V)
\simeq\det W\otimes\det V\simeq\det V\otimes\det W
\end{equation}
in which the first and third isomorphism are those given by
Lemma~\ref{lem:basicisom}
(with $V,W$ and $W,V$ respectively), $\det V\otimes\det W$ is mapped
into itself by multiplication by $(-1)^{\dim V\cdot\dim W}$.
\end{lemma}
\begin{proof}
This follows immediately from the familiar equation in anti-commuting
algebras (such as $\bigwedge^* V$),
\begin{equation}
a \wedge b = (-1)^{\deg a \cdot \deg b} b \wedge a
\end{equation}
which in turn follows from the equation for elements of degree 1
\begin{equation}
v \wedge w = - w \wedge v
\end{equation}
by induction.
\end{proof}

\section{Unordered Collections of Vector Spaces}

Now suppose we have an \emph{unordered} collection of $k$ vector spaces.
Formally, this is an element of the category $\Vect_k^{\text{i,u}}$
whose objects are
$k$-tuples of vector spaces $\{V_i\}_{i=1}^k$
and in which a morphism from $\{V_i\}_{i=1}^k$
to $\{W_i\}_{i=1}^k$ is a permuation $\pi:
\{1,\ldots,k\}\rightarrow\{1,\ldots,k\}$ and vector space
isomorphisms $l_i: V_i\rightarrow W_{\pi(i)}$.  There
is a functor from this category to $\Vect^{\text{i}}$, the direct sum, given by
\begin{equation}
\begin{split}
\{V_i\}_{i=1}^k \mapsto \bigoplus_{i=1}^k V_i \\
(\pi, \{l_i\}_{i=1}^k) \mapsto \sigma_\pi \circ \bigoplus_{i=1}^k l_i
\end{split}
\end{equation}
where $\sigma_\pi: \bigoplus_{i=1}^k W_{\pi(1)} \rightarrow \bigoplus_{i=1}^k
W_i$ is the canonical isomorphism rearranging the factors of
$\bigoplus_{i=1}^k W_i$.  Similarly, there is another functor to
$\Vect^{\text{i}}$, the tensor product, given by $\{V_i\}_{i=1}^k \mapsto
\bigotimes_{i=1}^k V_i$ and the corresponding map on morphisms.

Then, for any $\{V_i\}_{i=1}^k$ we can define the two 1-dimensional
vector spaces
\begin{equation}
\det (\bigoplus_{i=1}^k V_i)
\end{equation}
and
\begin{equation}
\bigotimes_{i=1}^k \det V_i
\end{equation}
by analogy with the construction above for $A \oplus B$.  But there is
no canonical isomorphism between them in general.  The story is more
complicated and more interesting.

\begin{lemma} \label{lem:evenrearr}
  If all the $\{V_i\}_{i=1}^k$ are even-dimensional, then there
is a canonical isomorphism 
\begin{equation}
\det (\bigoplus_{i=1}^k V_i) \simeq
\bigotimes_{i=1}^k \det V_i
\end{equation}
\end{lemma}
\begin{proof}
This is just an iterated version of Lemma~\ref{lem:basicisom}
above, with the further note that the isomorphism is independent of
the order of the $V_i$ by Lemma~\ref{lem:rearr}.
\end{proof}

For odd-dimensional vector spaces, we find

\begin{lemma} \label{lem:oddrearr}
If all the $\{V_i\}_{i=1}^k$ are odd-dimensional, there is a canonical
isomorphism
\begin{equation}
\det (\bigoplus_{i=1}^k V_i) \simeq
\bigotimes_{i=1}^k \det V_i \otimes \det(\bigoplus_{i=1}^k \RR v_i^*)
\end{equation}
where the $v_i$ are arbitrary independent vectors (so
$\bigoplus_{i=1}^k \RR v_i^*$ is
the dual of the free vector space generated by $\{1, \ldots,
k\}$)\footnote{Of course, this vector space has a canonical basis and
so is naturally isomorphic with its dual, but this formulation makes
the result easier to state.}.
\end{lemma}
\begin{proof}
Given an element $a_1 \otimes\cdots\otimes a_k \otimes\lambda$ of the
second vector space, the element
\begin{equation}
\lambda(v_1, \ldots, v_k)\cdot a_1 \wedge \cdots\wedge a_k
\end{equation}
(using the inclusions of $\det V_i$ in $\bigwedge^* (\bigoplus_{i=1}^k
V_i)$) is a non-zero element of $\det (\bigoplus_{i=1}^k V_i)$ which
pdoesn't depend on the particular ordering.  (By Lemma~\ref{lem:rearr},
both $\lambda(v_1, \ldots, v_k)$ and $a_1 \wedge\cdots\wedge a_k$
change sign under a simple transposition.)
\end{proof}

Since we only admit vector space isomorphisms, we can sort the vector
spaces by dimension.  In a somewhat weaker and more precise form, we
have a functor from 
$\Vect_k^{\text{i,u}}$ to $\Vect_2^{\text{i}}$ (the category of pairs of
vector spaces with isomorphisms) given by
\begin{equation} \label{eq:splitvect}
\bigoplus V_i \mapsto
\bigoplus_{\dim V_i \mathrm{\ even}}
\negthickspace\negthickspace V_i \oplus
\bigoplus_{\dim V_i \mathrm{\ odd}}
\negthickspace\negthickspace V_i
\end{equation}
(i.e., divide the vector spaces into even and odd dimensional pieces).

\begin{proposition} \label{lem:genrearr}
For arbitrary $V_i$, there is a canonical
isomorphism
\begin{equation}
\det (\bigoplus_{i=1}^k V_i) \simeq
\bigotimes_{i=1}^k \det V_i \otimes
\det\left(\bigoplus_{\dim V_i \mathrm{\ odd}}
\negthickspace\negthickspace\RR v_i^*\right)
\end{equation}
\end{proposition}
\begin{proof}
By equation~\ref{eq:splitvect}, Lemma~\ref{lem:basicisom}, and
Lemmas~\ref{lem:evenrearr} and~\ref{lem:oddrearr}, we have
\begin{equation}
\begin{split}
\det (\bigoplus_{i=1}^k V_i) &\simeq 
\det\left(\bigoplus_{\dim V_i \mathrm{\ even}}
\negthickspace\negthickspace V_i \oplus
\bigoplus_{\dim V_i \mathrm{\ odd}}
\negthickspace\negthickspace V_i\right)
\\
&\simeq
\det\left(\bigoplus_{\dim V_i \mathrm{\ even}}
\negthickspace\negthickspace V_i\right) \otimes
\det\left(\bigoplus_{\dim V_i \mathrm{\ odd}}
\negthickspace\negthickspace V_i\right)
\\
&\simeq
\bigotimes_{\dim V_i \mathrm{\ even}}
\negthickspace\negthickspace \det V_i \otimes
\bigotimes_{\dim V_i\,\mathrm{odd}}
\negthickspace\negthickspace \det V_i \otimes 
\det\left(\bigoplus_{\dim V_i \mathrm{\ odd}}
\negthickspace\negthickspace\RR v_i^*\right)
\\
&\simeq
\bigotimes_{i=1}^k \det V_i \otimes
\det\left(\bigoplus_{\dim V_i \mathrm{\ odd}}
\negthickspace\negthickspace\RR v_i^*\right)
\end{split}
\end{equation}
\end{proof}

\section {Exact Sequences}

In addition to the isomorphisms above, I will need some more general
isomorphisms involving $\det$ associated with any exact sequence.

\begin{lemma} \label{lem:detquot}
For any vector space $V$ of dimension $n$ and subspace $W$ of
dimension $m$, there is a canonical isomorphism
\begin{equation}
\det V \simeq \det W \otimes \det (V/W)
\end{equation}
\end{lemma}
\begin{proof}
As noted above, the natural injection $W \rightarrow V$ induces an
injection $\det W \rightarrow \bigwedge^m V$.
Similarly, the natural surjective map $V \rightarrow V/W$ induces a
surjective map $\bigwedge^{n-m} V \rightarrow \det (V/W)$,
whose kernel is generated by elements of the form $w \wedge v_1 \wedge
\cdots \wedge v_{n-m-1}$, $w \in W$.
For any two elements $a \in \det W$, $b \in \det (V/W)$, pick an
arbitrary preimage $b^* \in \bigwedge^{n-m}$ of $b$.  Then $a \wedge
b^*$ is independent of the choice of $b^*$, since $a \wedge c$ is 0
for any $c$ in the kernel of the map $\bigwedge^{n-m} V
\rightarrow \det (V/W)$ and defines a canonical (non-zero)
bilinear map from $\det W$ and $\det (V/W)$ to $\det V$, as desired.
\end{proof}

\begin{proposition} \label{prop:deteseq}
In any exact sequence
\begin{equation}
0 \stackrel{l_0}{\longrightarrow} V_1 \stackrel{l_1}{\longrightarrow}
\cdots \stackrel{l_{k-1}}{\longrightarrow} V_k
\stackrel{i_k}{\longrightarrow} 0
\end{equation}
of vector spaces, there is a canonical isomorphism
\begin{equation}
\bigotimes_
{\substack{
    1 \le i \le k \\
    k \mathrm{\ odd}
}}
\det V_i
\simeq
\bigotimes_
{\substack{
    1 \le i \le k \\
    k \mathrm{\ even}
}}
\det V_i
\end{equation}
\end{proposition}
\begin{proof}
For each $i$, we have, by Lemma~\ref{lem:detquot}, an isomorphism
\begin{equation}
\det V_i \simeq \det (\ker l_i) \otimes \det (V_i/\ker l_i)
\end{equation}
so, combining these (and using $\ker l_i \simeq V_{i-1}/\ker l_{i-1}$),
\begin{equation}
\begin{split}
\bigotimes_
{\substack{
    1 \le i \le k \\
    k \mathrm{\ odd}
}}
\det V_i
& \simeq
\bigotimes_
{\substack{
    1 \le i \le k \\
    k \mathrm{\ odd}
}}
\det (\ker l_i) \otimes \det (V_i/\ker l_i) \\
& \simeq
\bigotimes_
{\substack{
    1 \le i \le k \\
    k \mathrm{\ odd}
}}
\det (V_{i-1}/\ker l_{i-1}) \otimes \det \ker l_{i+1} \\
& \simeq
\bigotimes_
{\substack{
    1 \le i \le k \\
    k \mathrm{\ even}
}}
\det V_i
\end{split}
\end{equation}
\end{proof}

%% file: graph.tex
\chapter{Graphs}
\label{sec:graph}

In this Appendix I formally prove the equivalence of three different
notions of orientation on a graph, generalizing the comment in
Section~\ref{sec:integrals} about the two different ways of writing
the integral invariants.

An (unoriented) \emph{graph based on a circle} is an undirected graph
$\Gamma$
with a choice of a (directed) cycle $c(\Gamma)$ with distinct edges
and vertices.
An \emph{orientation} on a graph based on a circle $\Gamma$ is an
orientation on the 
1-dimensional vector space
\begin{equation}
\orient(\Gamma) = \det \RR V(\Gamma) \otimes \bigotimes_{e \in
E(\Gamma)} \RR V(e)
\end{equation}
where
\begin{itemize}
\item $\det W$, for $W$ an $n$-dimensional vector space, is
$\bigwedge^n W$.
\item $V(\Gamma)$ is the set of vertices of $\Gamma$.
\item $E(\Gamma)$ is the set of edges of $\Gamma$.
\item $VE(\Gamma)$ is the set of incident vertex-edge pairs of $\Gamma$.
\item $V(e)$ is the set of vertices that are endpoints of $e$ (a two
element set).
\item $E(v)$, which we will need later, is the set of edges that
have the vertex $v$ as an endpoint.
\item $V(c)$ and $E(c)$ are the vertices and edges, respectively, in
$c(\Gamma)$.
\item $V_i(\Gamma)$ and $E_i(\Gamma)$, the internal vertices and edges
of $\Gamma$, are $V(\Gamma) \\ V(c)$ and  $E(\Gamma) \\ E(c)$,
respectively.
\end{itemize}
Two equivalent definitions are given in Proposition~\ref{prop:or}
below.

A \emph{labelling} $l$ of the
internal edges of $\Gamma$ is a bijective map
\begin{equation}
l: \{1, \ldots, |E_i(\Gamma)|\} \rightarrow E_i(\Gamma)
\end{equation}
together with bijective maps
\begin{equation}
l_e: \{s, d\} \rightarrow V(e), \quad e\in E_i(\Gamma)
\end{equation}

A labelling $l$ and an orientation $o$ together define a chain in
$(S^2)^{E_i(\Gamma)}$ for any given knot $K$.  Pick an arbitrary bijection
\begin{equation}
t: V(\Gamma) \rightarrow \{1, \ldots, |V(\Gamma)|\}
\end{equation}
Then the chain in $(S^2)^{E_i(\Gamma)}$ is
\begin{multline}
\phi_{\Gamma,o,l}(K) = \epsilon
	\phi_{l_{l(1)}(s),l_{l(1)}(d)} (C_{|V(c)|,|V_i(\Gamma)|}^t(K)) \times\\
	\phi_{l_{l(2)}(s),l_{l(2)}(d)} (C_{|V(c)|,|V_i(\Gamma)|}^t(K)) \times
	\cdots \times
	\phi_{l_{l(|E_i(\Gamma)|)}(s),l_{l(|E_i(\Gamma)|)}(d)}
				(C_{|V(c)|,|V_i(\Gamma)|}^t(K))
\end{multline}
where 
\begin{itemize}
\item $C_{|V(c)|,|V_i(\Gamma)|}^t(K)$ is the component of
$C_{|V(c)|,|V_i(\Gamma)|}(K)$ with corresponding labellings of the points
and with the points on the circle occuring in the same cyclic order as
they do in $\Gamma$, imbued with the orientation the interior inherits
as a subset of $(S^1)^{|V(c)|} \times (\RR^3)^{|V_i(\Gamma)|}$ with
the factors in the appropriate order.
\item $\epsilon$ is the appropriate sign: $\pm 1$, depending on whether
or not the element
\begin{equation}
(t(1) \wedge \cdots \wedge t(n)) \otimes
	\bigotimes_{e\in V_i(\Gamma)} (l_e(1) \wedge l_e(2))
	\otimes \bigotimes_{e\in V(c)} (c_e(1) \wedge c_e(2))
	\in \orient(\Gamma)
\end{equation}
is positively oriented.  ($c_e: \{1,2\}\rightarrow V(e)$ is the map
determining the orientation of the circle.)
\end{itemize}

Note, in the configuration space of the knot, that cycle vertices
are constrained to lie on the knot, internal vertices are not so
constrained, internal edges give rise to propogators, while cylce
edges are merely place-holders to indicate the order of points around
the knot.

This is a good definition: it doesn't depend on the labelling $t$ of
$V(\Gamma)$ we chose, because of the choice of $\epsilon$.  (Note in
particular that both $S^1$ and $\RR^3$ are odd-dimensional, so the
orientation of $C_{|V(c)|,|V_i(\Gamma)|}^t$ depends on the sign of $t$
(considered as a permutation),
exactly as $\epsilon$ does.)

In addition, if given an orientation $o$ and a labelling $t$ of the
vertices as above, there is a well-defined form
\begin{equation}
\theta_{\Gamma,o,t} = \phi_{\Gamma,o,l}^* \omega^{\times|E_i(\Gamma)|}
\end{equation}
on $C_{|V(c)|,|V_i(\Gamma)|}^t$
independent of the labelling of the edges $l$.  With the orientation
defined by $t$,
the integral of this form is the volume of the chain
$\phi_{\Gamma,o,l}$ and is thus independent of $t$ and $l$.  These
integrals are exactly those defined in BT.

These integrals and chains are the basic building blocks to construct
knot invariants.  But there are two important equivalent ways of
defining the orientation above.

\begin{proposition} \label{prop:or}
Orientations on the following 1-dimensional vector spaces are
equivalent:
\begin{enumerate}
\item $\orient(\Gamma) = \det \RR V(\Gamma) \otimes
	\bigotimes_{e \in E(\Gamma)} \RR V(e)$
\item $\bigotimes_{v \in V(\Gamma)} \det \RR E(v) \otimes
		\det(\oplus_{\deg v~\text{even}} \RR v)$
\item $\det\RR E(\Gamma) \otimes \det H_1(\Gamma) \otimes \det H_0(\Gamma)$
\end{enumerate}
\end{proposition}

Note that we have the following equivalences:
\begin{equation} \label {eq:or1}
\begin{aligned}
\bigotimes_{e \in E(\Gamma)} \det \RR V(e) &\simeq
	\det \RR VE(\Gamma)\\
	&\simeq 
	\bigotimes_{v \in V(\Gamma)} \det \RR E(v) \otimes
		\det (\bigoplus_{\deg v~\text{odd}} \RR v)
\end{aligned}
\end{equation}
The first isomorphism arises from grouping $VE(\Gamma)$ by vertex-edge
pairs by those sharing an edge, while the second comes about by
grouping by vertices.  Note that when we group vertex-edge pairs in
this way, we must be careful to observe the ordering of the groups.  A
group with an even number of elements will commute with all other
groups, while a group with an odd number of elements will
anti-commute.  Every edge has two endpoints, yielding the first
isomorphism, while the second isomorphism is just the above
observation.  Tensoring equation~\ref{eq:or1} with $\det V(\Gamma)$
yields $(1) \simeq (2)$.

To see $(1) \simeq (3)$, note that, since a graph has no two-cells, we
have an exact sequence
\begin{equation} \label{eq:exactseq}
0 \rightarrow H_1(\Gamma) \rightarrow C_1(\Gamma) \rightarrow C_0(\Gamma)
	\rightarrow H_0(\Gamma) \rightarrow 0
\end{equation}
yielding an isomorphism
\begin{equation} \label{eq:or2}
\det H_1(\Gamma) \otimes \det H_0(\Gamma) \simeq
	\det C_1(\Gamma) \otimes \det C_0(\Gamma)
\end{equation}
(all these spaces have canonical (unordered) bases, so there's no
worry about dual spaces).  $C_0(\Gamma)$ is just $\RR V(\Gamma)$, but
$C_1(\Gamma)$ is the vector space spanned by \emph{oriented} edges, so
\begin{gather}
C_1(\Gamma) \simeq \bigoplus_{e \in E(\Gamma)} \det V(e)\\
\det C_1(\Gamma) \simeq \det \RR E(\Gamma) \otimes
	\bigotimes_{e \in E(\Gamma)} \det V(e)
\end{gather}
Therefore, tensoring equation~\ref{eq:or2} above with $\det \RR
E(\Gamma)$, we find
\begin{equation}
\det\RR E(\Gamma) \otimes \det H_1(\Gamma) \otimes \det H_0(\Gamma)
	\simeq \det \RR V(\Gamma) \otimes \bigotimes_{e \in E(\Gamma)} \RR V(e)
\end{equation}
as desired.

Note that orientations of type (2) in the proposition correspond to
Bar-Natan's trivalent diagrams, as in BN.  Orientations of type (3)
are mentioned in Kontsevich's paper K.  (He omits the $H_0$ term.  In
the case of knots, this is OK, as all integrals of interest are
connected (to the knot).  However, to get invariants of homology
3-spheres the $H_0$ term is necessary.)

Also note that in each of the cases above the cycle of the boundary
could have been omitted from the graph $\Gamma$ without essentially
changing the definition.  For instance, for definition (1) of
$\orient(\Gamma)$, if we omit the circle from the graph we omit no
vertices and some edges, each with a canonical orientation.  This
would have made the definition of the chain $\phi_{\Gamma,o,l}$
slightly cleaner.  However,
the definition of the boundary operator is easier to state if the
cycle is included in the graph.